\documentclass[a4paper]{article}%
\usepackage{amsmath,amsthm,amssymb,amsfonts}
\usepackage{enumerate,amscd,amsxtra}
\usepackage{mathrsfs}
\usepackage{graphicx}
\usepackage[cmtip,all]{xy}
\usepackage{amsmath}
\usepackage{amsfonts}
\usepackage{amssymb}%
\providecommand{\U}[1]{\protect\rule{.1in}{.1in}}
\providecommand{\U}[1]{\protect\rule{.1in}{.1in}}
\providecommand{\U}[1]{\protect\rule{.1in}{.1in}}
\normalsize
\newtheoremstyle{theoremstyle}
{10pt}
{5pt}
{\itshape}
{}
{\bfseries}
{:}
{.5em}
{}
\newtheoremstyle{examplestyle}
{10pt}
{5pt}
{}
{}
{\bfseries}
{:}
{.5em}
{}
\theoremstyle{theoremstyle}
\newtheorem{theorem}{Theorem}[section]
\newtheorem*{theorem*}{Theorem}
\newtheorem{lemma}[theorem]{Lemma}
\newtheorem{proposition}[theorem]{Proposition}
\newtheorem*{proposition*}{Proposition}
\newtheorem{corollary}[theorem]{Corollary}
\newtheorem*{corollary*}{Corollary}

\theoremstyle{definition}

\newtheorem{examples}[theorem]{Examples}

\newtheorem{definition*}{Definition}
\newtheorem{remark}[theorem]{Remark}
\newtheorem{remark*}{Remark}

\begin{document}

\title{\textbf{Hermitian $K$-theory and $2$-regularity for totally real number
fields}}
\author{A.\thinspace J.~Berrick, M.~Karoubi, P.\thinspace A.~{\O }{}stv{\ae }{}r}
\date{\today}
\maketitle

\begin{abstract}
We completely determine the $2$-primary torsion subgroups of the hermitian
$K$-groups of rings of $2$-integers in totally real $2$-regular number fields.
The result is almost periodic with period $8$. Moreover, the $2$-regular case
is precisely the class of totally real number fields that have homotopy
cartesian \textquotedblleft B\"{o}kstedt square", relating the $K$-theory of
the $2$-integers to that of the fields of real and complex numbers and finite
fields. We also identify the homotopy fibers of the forgetful and hyperbolic
maps relating hermitian and algebraic $K$-theory. The result is then exactly
periodic of period $8$ in the orthogonal case. In both the orthogonal and
symplectic cases, we prove a $2$-primary hermitian homotopy limit conjecture
for these rings.

\end{abstract}

\newpage

\section{Introduction and statement of results}

\label{section:mainresults}

Let $F$ be a real number field with $r$ real embeddings, ring of integers
$\mathcal{O}_{F}$ and ring of $2$-integers $R_{F}=\mathcal{O}_{F}[1/2]$. A key
ingredient in the study of the $K$-theory of these rings is the commuting
square of $2$-completed connective\textsf{ }$K$-theory spectra%
\begin{equation}%
\begin{array}
[c]{ccc}%
\mathcal{K}(R_{F})_{\#} & \rightarrow & {\displaystyle\bigvee\limits^{r}}%
{}\mathcal{K}(\mathbb{R})_{\#}^{c}\\
\downarrow &  & \downarrow\\
\mathcal{K}(\mathbb{F}_{q})_{\#} & \rightarrow & {\displaystyle\bigvee
\limits^{r}}\mathcal{K}(\mathbb{C})_{\#}^{c}%
\end{array}
\label{algebraic K Bokstedt square}%
\end{equation}
introduced by B\"{o}kstedt for the rational case ($r=1$) in \cite{Bok}, and in
\cite{HO}, \cite{Mitchell} for the general case; see also Appendix A. Here,
$\mathbb{F}_{q}$ is a finite residue field described later in this
Introduction. The $K$-theory spectrum $\mathcal{K}(\Lambda)$ refers to the
spectrum defined by $\mathcal{K}(\Lambda)_{n}=\Omega{B\mathrm{GL}(S}%
^{n+1}{\Lambda)}^{+}$ for $n\geq0$, where $S^{m}\Lambda$ denotes the
$m$-iterated suspension of $\Lambda$ (when $\Lambda$ is a discrete ring) or
Calkin algebra ($\Lambda$ a topological ring), as defined in \cite[Appendix
A]{BK}. Throughout this paper, we use the notation \thinspace$\mathcal{E}^{c}$
to denote the connective covering of a spectrum $\mathcal{E}$, and $_{\#}$ to
denote $2$-adic completion of connective spectra or groups. The often-used
notation $\mathcal{E}_{\#}^{c}$ means $(\mathcal{E}^{c})_{\#}.$ For instance,
\[
\mathcal{K}(\mathbb{R})_{\#}^{c}=(\mathcal{K}(\mathbb{R})^{c})_{\#}%
\quad\text{and\quad}\mathcal{K}(R_{F})_{\#}^{c}=(\mathcal{K}(R_{F})^{c}%
)_{\#}=\mathcal{K}(R_{F})_{\#}%
\]
(see also the discussion after Lemma \ref{Connective.Spectra}). The bottom
horizontal map is the Brauer lift, corresponding to the fibring of Adams' map
$\psi^{q}-1$ on the $2$-completed connective complex topological $K$-theory
spectrum $\mathcal{K}(\mathbb{C})_{\#}^{c}$. The remaining maps are induced
from the obvious ring homomorphisms via Suslin's identification of the
$2$-completed algebraic $K$-theory spectra of the real and complex numbers
with $\mathcal{K}(\mathbb{R})_{\#}^{c}$ and $\mathcal{K}(\mathbb{C})_{\#}^{c}%
$, respectively, in \cite{Suslin :local fields}. If preferred, one can think
of the above in terms of spaces and maps; however, the spaces also have an
infinite loop space structure that is preserved by the maps.

In the rational case, the Dwyer-Friedlander formulation of the
Quillen-Lichtenbaum conjecture for $\mathbb{Z}$ at the prime $2$ is that the
above square is homotopy cartesian, see \cite[Conjecture 1.3, Proposition
4.2]{DF}. This has been affirmed in work of Rognes and Weibel \cite{RW} (see
\cite[Corollary 8]{W:CR}), as a consequence of \cite{Bok}, Voevodsky's
solution of the Milnor Conjecture \cite{V} and his subsequent joint work with
Suslin \cite{SV}.

In the general number field case, Rognes and Weibel \cite{RW} determined the
groups $K_{n}(R_{F})_{\#}$ up to extensions. It turns out that in the case of
$2$-regular real number rings, discussed below, these extension problems
disappear \cite{RO}. This leads to the above square being homotopy cartesian
in that case too. Many of the foregoing results having been developed for
spaces, we present the transition to spectra in Appendix A.
These developments raise the question of for which class of real number fields
$F$ the square (\ref{algebraic K Bokstedt square}) is homotopy cartesian.

We turn now to the hermitian analog of the above. For the definition of
hermitian $K$-theory we refer to \cite{K:AnnM112hgo} and \cite[Introduction
and Appendix A]{BK}. Briefly, for a ring $\Lambda$ with involution
${}_{\varepsilon}{KQ}_{0}{(\Lambda)}$ denotes the Grothendieck group of
isomorphism classes of finitely generated projective $\Lambda$-modules with
nondegenerate $\varepsilon$-hermitian form, where we let $\varepsilon=\pm1$
according to whether orthogonal ($\varepsilon=+1$) or symplectic
($\varepsilon=-1$) actions on the ring $\Lambda$ are involved. If $\Lambda$ is
discrete, $B{{}_{\varepsilon}O(\Lambda)}^{+}$ represents the plus-construction
of the classifying space of the limit ${{}_{\varepsilon}O(\Lambda)}$ of the
$\varepsilon$\emph{-}orthogonal\emph{ }groups ${}_{\varepsilon}O_{n,n}%
(\Lambda)$. This last group is the group of automorphisms of the $\varepsilon
$-hyperbolic module ${}_{\varepsilon}H(\Lambda^{n})$, whose elements can be
described as $2\times2$ matrices written in $n$-blocks%
\[
M=\left[
\begin{array}
[c]{cc}%
a & b\\
c & d
\end{array}
\right]
\]
such that $M^{\ast}M=MM^{\ast}=I$, where the \textquotedblleft$\varepsilon
$-hyperbolic adjoint\textquotedblright\ $M^{\ast}$ is defined as
\[
M^{\ast}=\left[
\begin{array}
[c]{cc}%
^{\mathrm{t}}\check{d} & \varepsilon\,^{\mathrm{t}}\check{b}\\
\check{\varepsilon}\,^{\mathrm{t}}\check{c} & ^{\mathrm{t}}\check{a}%
\end{array}
\right]  \text{.}%
\]
For a discrete ring $A$, the $\varepsilon$-hermitian $K$-theory spectrum
${}_{\varepsilon}\mathcal{KQ}(A)$ refers to the spectrum defined by
${}_{\varepsilon}\mathcal{KQ}(A)_{n}=\Omega{B{}_{\varepsilon}O(S}^{n+1}%
{A)}^{+}$ for $n\geq0$, where $S^{m}A$ denotes the $m$-iterated suspension of
$A$, as defined in \cite[Appendix A]{BK}. We note that ${}_{\varepsilon
}\mathcal{KQ}(A)_{0}$ has non-naturally the homotopy type of ${}_{\varepsilon
}KQ_{0}(A)\times B{}_{\varepsilon}O(A)^{+}$, where ${}_{\varepsilon}KQ_{0}(A)$
is endowed with the discrete topology. The same definition applies for the
$K$-theory spectrum $\mathcal{K}(A)$ on replacing the orthogonal group by the
general linear group. There is however a significant difference between the
two theories, at least for regular noetherian rings like fields or Dedekind
rings. For such rings $A$, $\mathcal{K}(A)=\mathcal{K}(A)^{c}$ is connective
\textsl{i.e.} $K_{n}(A)=0$ for $n<0$, whereas ${}_{\varepsilon}\mathcal{KQ}%
(A)$ is not connective in general.

On the other hand, the spectra of topological hermitian $K$-theory (with
trivial involutions on $\mathbb{R}$ and $\mathbb{C}$) ${}_{\varepsilon
}\mathcal{KQ}(\mathbb{R})$ and ${}_{\varepsilon}\mathcal{KQ}(\mathbb{C})$ have
been defined in \cite[Appendix A]{BK}. A geometric description of these
spectra appears in Appendix B below. In order to avoid a potential confusion
between hermitian $K$-theory and surgery theory, we are writing ${}%
_{\varepsilon}\mathcal{KQ}$\ for the hermitian $K$-theory spectrum, and
${}_{\varepsilon}{KQ}_{n}$\ for the corresponding homotopy groups. (These are
denoted by ${}_{\varepsilon}\mathcal{L}$\ and ${}_{\varepsilon}{L}_{n}%
$\ respectively in \cite{BK}.)

In \cite{BK}, the first two authors constructed a Brauer lift in hermitian
$K$-theory and considered the hermitian analogue of the B\"{o}kstedt square
for the rational numbers $\mathbb{Q}$, \textsl{i.e.} for $r=1$ and
$R_{F}=\mathbb{Z}[1/2]$. For general number fields, the commuting B\"{o}kstedt
square for hermitian $K$-theory takes the form:
\begin{equation}%
\begin{array}
[c]{ccc}%
{}_{\varepsilon}\mathcal{KQ}(R_{F})_{\#}^{c} & \rightarrow &
{\displaystyle \bigvee\limits^{r}}{}_{\varepsilon}\mathcal{KQ}(\mathbb{R}%
)_{\#}^{c}\\
\downarrow &  & \downarrow\\
{}_{\varepsilon}\mathcal{KQ}(\mathbb{F}_{q})_{\#}^{c} & \rightarrow &
{\displaystyle\bigvee\limits^{r}}{}_{\varepsilon}\mathcal{KQ}(\mathbb{C}%
)_{\#}^{c}%
\end{array}
\label{hermitian K Bokstedt square}%
\end{equation}

It was shown in \cite{BK} that when $F=\mathbb{Q}$ the square
(\ref{hermitian K Bokstedt square}) too is homotopy cartesian; the result
leads to another version of the homotopy limit problem related to the
Quillen-Lichtenbaum conjecture, expressed as the $2$-adic homotopy equivalence
of \emph{ }the fixed point set and the homotopy fixed point set of the
${}_{\varepsilon}\mathbb{Z}/2$ action on $\mathcal{K}(\mathbb{Z}[1/2])$. Thus
again, one is led to ask \textit{for which class of totally real number fields
the square (\ref{hermitian K Bokstedt square}) is homotopy cartesian, and for
which the homotopy equivalence generalizes.} \smallskip

These are the principal questions addressed in the present work. \smallskip

Tackling these questions leads to a focus on a particular class of number
fields $F$ with the associated rings of integers $\mathcal{O}_{F}$, as
follows. From a theorem of Tate \cite[Theorem 6.2]{Tate}, one knows that the
$2$-primary part of the finite abelian group $K_{2}(\mathcal{O}_{F})$ has
order at least $2^{r}$, where $r$ is the number of real embeddings. We call
$F$ (and $\mathcal{O}_{F}$, $R_{F}$) $2$\emph{-regular }when this order is
exactly $2^{r}$. See Proposition \ref{2+-regular characterization} below for
alternative characterizations. In the totally real case, which is our concern
here, $r=[F:\mathbb{Q}]$. The simplest examples are the rational numbers
$\mathbb{Q}$ and the following fields recorded in \cite[\S 4]{RO}.

\begin{enumerate}
\item Let $b\geq2$. The maximal real subfield $F=\mathbb{Q}(\zeta_{2^{b}}%
+\bar{\zeta}_{2^{b}})$ of $\mathbb{Q}(\zeta_{2^{b}})$ is a totally real
$2$-regular number field with $r=2^{b-2}$.

\item Let $m$ be an odd prime power such that $2$ is a primitive root modulo
$m$. Then $F=\mathbb{Q}(\zeta_{m}+\bar{\zeta}_{m})$ is a totally real
$2$-regular number field when Euler's $\phi$-function $\phi(m)\leq66$ (except
for $m=29$), and also for Sophie Germain primes ($m$ and $(m-1)/2$ both prime)
with $m\not \equiv 7\;(\mathrm{mod}\ 8)$ (the first few instances are
$m=5,11,59,83,107$ and $179$). The number $r$ of real embeddings is
$\phi(m)/2$.

\item Let $F=\mathbb{Q}(\sqrt{d})$ be a quadratic number field with $d>0$
square free. Then $F$ is $2$-regular if and only if $d=2$, $d=p$ or $d=2p$
with $p\equiv\pm3\;(\mathrm{mod}\ 8)$ prime \cite{BS}. Here, $r=2$.
\end{enumerate}

\smallskip

The residue field $\mathbb{F}_{q}$ of $R_{F}$ referred to above is now chosen
in the following manner. The number $q$ is a prime number with this property:
the elements corresponding to the Adams operations $\psi^{q}$ and $\psi^{-1}$
in the ring of operations of the periodic complex topological $K$-theory
spectrum generate the Galois group of $F(\mu_{2^{\infty}}(\mathbb{C}))/F$
obtained by adjoining all $2$-primary roots of unity $\mu_{2^{\infty}%
}(\mathbb{C})\subset\mathbb{C}$ to $F$ \cite[\S 1]{Mitchell}. The Cebotarev
density theorem guarantees the existence of infinitely many such prime powers.
By Dirichlet's theorem on arithmetic progressions we may assume that $q$ is a
prime number, an hypothesis we assume throughout all the paper.
According to \cite{Mitchell}, if $a_{F}:=(|\mu_{2^{\infty}}(F(\sqrt
{-1}))|)_{2}$ is the $2$-adic valuation, then $q$ is $\equiv\pm
1\;(\mathrm{mod}\ 2^{a})$ but not $(\mathrm{mod}\ 2^{a+1})$. In the examples
above: when $F=\mathbb{Q}(\zeta_{2^{b}}+\bar{\zeta}_{2^{b}})$, then $a_{F}=b$;
when $F=\mathbb{Q}(\zeta_{m}+\bar{\zeta}_{m})$ or $\mathbb{Q}(\sqrt{d})$ with
$d>2$, $a_{F}=2$; and finally, when $F=\mathbb{Q}(\sqrt{2})=\mathbb{Q}%
(\zeta_{8}+\bar{\zeta}_{8})$, we have $a_{F}=3$.

\smallskip

We are now ready to state our main results.

\begin{theorem}
\label{theorem1} For every totally real $2$-regular number field $F$, and for
any $q$ as discussed above, the square (\ref{hermitian K Bokstedt square}) is
homotopy cartesian for $\varepsilon=\pm1$.
\end{theorem}

Our proof of this theorem is based on the techniques employed in the case of
the rational numbers \cite{BK} and the analogous algebraic $K$-theoretic
result established in \cite{HO}, \cite{Mitchell} and \cite{RO} (see Appendix A
for an overview). \vspace{0.1in}

The next result crystallizes the special role of $2$-regular fields in this setting.

\begin{theorem}
\label{converse to theorem 1}Let $q$ be as above. Then, for every totally real
number field $F$, the following are equivalent.

\begin{enumerate}
\item[(i)] $F$ is $2$-regular.

\item[(ii)] The square (\ref{hermitian K Bokstedt square}) is homotopy
cartesian for $F$ when $\varepsilon=1$.

\item[(iii)] The square (\ref{algebraic K Bokstedt square}) is homotopy
cartesian for $F$.
\end{enumerate}
\end{theorem}

Since the Quillen-Lichtenbaum conjecture has been established for every real
number field by the third author \cite{Ostvar}, one consequence of this
theorem is that, in contrast to the rational case, for general real number
fields the Quillen-Lichtenbaum conjecture fails to imply that the square
(\ref{algebraic K Bokstedt square}) is homotopy cartesian.

\smallskip

Theorem \ref{theorem1} allows us to compute explicitly the $2$-primary torsion
of the nonnegative hermitian $K$-groups ${}_{\varepsilon}KQ_{n}(R_{F})$ for
$F$ as above. We tabulate and compare these groups with the corresponding
algebraic $K$-groups $K_{n}(R_{F})$ computed in \cite{HO} and \cite{RO}.

\begin{theorem}
\label{theorem2}Let $F$ be a totally real $2$-regular number field. Up to
finite groups of odd order, the groups ${}_{\varepsilon}KQ_{n}(R_{F})$ are
given in the following table. (If $m$ is even, let $w_{m}=2^{a_{F}+\nu_{2}%
(m)}$; also, $\delta_{n0}$ denotes the Kronecker symbol.) \begin{table}[tbh]
\begin{center}%
\begin{tabular}
[c]{p{0.4in}|p{1.5in}|p{1.1in}|p{1.5in}|}\hline
$n\geq0$ & ${}_{-1}KQ_{n}(R_{F})$ & ${}_{1}KQ_{n}(R_{F})$ & $K_{n}(R_{F}%
)$\\\hline
$8k$ & $\delta_{n0}\mathbb{Z}$ & $\delta_{n0}\mathbb{Z}\oplus\mathbb{Z}%
^{r}\oplus\mathbb{Z}/2$ & $\delta_{n0}\mathbb{Z}$\\
$8k+1$ & $0$ & $(\mathbb{Z}/2)^{r+2}$ & $\mathbb{Z}^{r}\oplus\mathbb{Z}/2$\\
$8k+2$ & $\mathbb{Z}^{r}$ & $(\mathbb{Z}/2)^{r+1}$ & $(\mathbb{Z}/2)^{r}$\\
$8k+3$ & $(\mathbb{Z}/2)^{r-1}\oplus\mathbb{Z}/2w_{4k+2}$ & $\mathbb{Z}%
/w_{4k+2}$ & $(\mathbb{Z}/2)^{r-1}\oplus\mathbb{Z}/2w_{4k+2}$\\
$8k+4$ & $(\mathbb{Z}/2)^{r}$ & $\mathbb{Z}^{r}$ & $0$\\
$8k+5$ & $\mathbb{Z}/2$ & $0$ & $\mathbb{Z}^{r}$\\
$8k+6$ & $\mathbb{Z}^{r}$ & $0$ & $0$\\
$8k+7$ & $\mathbb{Z}/w_{4k+4}$ & $\mathbb{Z}/w_{4k+4}$ & $\mathbb{Z}/w_{4k+4}%
$\\\hline
\end{tabular}
\end{center}
\end{table}
\end{theorem}

The proof of Theorem \ref{theorem2} makes use of a splitting result for
${}_{\varepsilon}\mathcal{KQ}(R_{F})^{c}$ shown in
\S \ref{section:splittingresults} and an explicit computation carried out in
\S \ref{section:proofoftheorem2}.

We recall that the forgetful and hyperbolic functors induce the two homotopy
fiber sequences
\[
{}_{\varepsilon}\mathcal{V}(R_{F})\longrightarrow{}_{\varepsilon}%
\mathcal{KQ}(R_{F})\longrightarrow\mathcal{K}(R_{F})\text{ and }%
{}_{\varepsilon}\mathcal{U}(R_{F})\longrightarrow\mathcal{K}(R_{F}%
)\longrightarrow{}_{\varepsilon}\mathcal{KQ}(R_{F})\text{.}%
\]
The fundamental theorem in hermitian $K$-theory \cite{K:AnnM112fun} states
that there is a natural homotopy equivalence
\[
{}_{\varepsilon}\mathcal{V}(R_{F})\simeq\Omega{}_{-\varepsilon}\mathcal{U}%
(R_{F})\text{.}%
\]

Our next result gives an explicit computation of the homotopy groups of these spectra.

\begin{theorem}
\label{theorem4} For any totally real $2$-regular number field $F$, the
groups
\[
{}_{\varepsilon}V_{n}(R_{F}):=\pi_{n}({}_{\varepsilon}\mathcal{V}(R_{F}%
))\cong\pi_{n}(\Omega{}_{-\varepsilon}\mathcal{U}(R_{F}))=:{}_{-\varepsilon
}U_{n+1}(R_{F})
\]
are given, up to finite groups of odd order, as below.\begin{table}[tbh]
\begin{center}%
\begin{tabular}
[c]{p{0.4in}|p{1.0in}|p{1.2in}|}\hline
$n\geq0$ & ${}_{-1}V_{n}(R_{F})$ & ${}_{1}V_{n}(R_{F})$\\\hline
$8k$ & $\mathbb{Z}^{r}\oplus\mathbb{Z}/2$ & $\mathbb{Z}^{2r}$\\
$8k+1$ & $0$ & $(\mathbb{Z}/2)^{2r}$\\
$8k+2$ & $\mathbb{Z}^{r}$ & $(\mathbb{Z}/2)^{2r}$\\
$8k+3$ & $0$ & $0$\\
$8k+4$ & $\mathbb{Z}^{r}$ & $\mathbb{Z}^{2r}$\\
$8k+5$ & $\mathbb{Z}/2$ & $0$\\
$8k+6$ & $\mathbb{Z}^{r}\oplus\mathbb{Z}/2$ & $0$\\
$8k+7$ & $\mathbb{Z}/2$ & $0$\\\hline
\end{tabular}
\end{center}
\end{table}

More precisely, the spectrum ${}_{1}\mathcal{V}(R_{F})_{\#}^{c}$ has the
homotopy type of
\[
\bigvee^{2r}\mathcal{K}(\mathbb{R})_{\#}^{c}\simeq{}_{1}\mathcal{V}%
(R_{\mathbb{Q}})_{\#}^{c}\vee\bigvee^{2(r-1)}\mathcal{K}(\mathbb{R})_{\#}%
^{c}\,\text{,}%
\]
The cup-product with a generator of $K_{8}(\mathbb{R})$ induces a periodicity
homotopy equivalence%
\[
{}_{1}\mathcal{V}(R_{F})_{\#}^{c}\simeq(\Omega^{8}{}_{1}\mathcal{V}%
(R_{F}))_{\#}^{c}\,\text{.}%
\]

\end{theorem}

It would be interesting to have more information about the class of number
fields for which the intriguing periodicity result in this theorem holds. Note
also that for $\varepsilon=-1$, the homotopy type of ${}_{\varepsilon
}\mathcal{V}(R_{F})_{\#}^{c}$ is more complicated to determine explicitly,
although we know that its homotopy groups are periodic (in Section
\ref{V-computation} and Appendix D we show that an analogous splitting of the
spectrum ${}_{-1}\mathcal{V}(R_{F})_{\#}^{c}$ as the product of classical
topological spectra does not hold).

Recall from \cite[\S 7]{BK} that the standard ${}_{\varepsilon}\mathbb{Z}%
/2$-action on $\mathcal{K}(R_{F})$ defined via conjugation of the involute
transpose of a matrix by
\[
{}_{\varepsilon}J_{n}=\left[
\begin{array}
[c]{cc}%
0 & \varepsilon I_{n}\\
I_{n} & 0
\end{array}
\right]
\]
induces an isomorphism between the fixed point spectrum $\mathcal{K}%
(R_{F})^{{}_{\varepsilon}\mathbb{Z}/2}$ and ${}_{\varepsilon}\mathcal{KQ}%
(R_{F})$.

Our next result solves in the affirmative a $2$-primary homotopy limit problem
for ${}_{\varepsilon}\mathbb{Z}/2$-homotopy fixed point $K$-theory spectra, in
the special case of totally real $2$-regular number fields. (In \cite{Ostvar}
the corresponding problem in algebraic $K$-theory was solved in the
affirmative for every real number field.) Recall the homotopy fixed point
spectrum
\[
\mathcal{K}(R_{F}){}^{h({}_{\varepsilon}\mathbb{Z}/2)}:=\mathrm{map}%
_{{}_{\varepsilon}\mathbb{Z}/2}(\Sigma^{\infty}E(\mathbb{Z}/2)_{+}%
,\,\mathcal{K}(R_{F}))\text{.}%
\]
Here $\mathrm{map}_{{}_{\varepsilon}\mathbb{Z}/2}$ denotes the function
spectrum of ${}_{\varepsilon}\mathbb{Z}/2$-equivariant maps and $E(\mathbb{Z}%
/2)$ is a free contractible $\mathbb{Z}/2$-space, such as the CW-complex
$S^{\infty}$ with antipodal action. Our hermitian analogue is now the following:

\begin{theorem}
\label{theorem3} For every totally real $2$-regular number field $F$ there is
a natural homotopy equivalence of $2$-completed spectra
\[
{}_{\varepsilon}\mathcal{KQ}(R_{F})_{\#}^{c}\simeq(\mathcal{K}(R_{F}){}%
_{\#}^{h({}_{\varepsilon}\mathbb{Z}/2)})^{c}.
\]

\end{theorem}

The proof of Theorem \ref{theorem3} follows from Theorem \ref{theorem1} and
results established in \cite{BK}. Earlier results in this direction inspired a
much more general conjecture formulated by B.~Williams in \cite[3.4.2]%
{Williams}. However, in Appendix C below, we provide a counterexample to that
conjecture, in the form of a ring with vanishing $K$-theory but nontrivial
$KQ$-theory.\smallskip

\begin{remark}
Most of the theorems in this introduction are also true if we replace the
$2$-completions of the spectra involved by their $2$-localizations. As we
shall see through the paper, the proofs of most of our main theorems work as
well in this context. There is an important exception however, namely our
Theorem \ref{theorem3}, which is only true for completions.\vspace{0.1in}
\end{remark}

In our work in progress \cite{BKO:bottperiodicity} the results shown in this
paper are used to give the first algebraic examples of Bott periodicity
isomorphisms for hermitian $K$-groups. Another project uses these results to
provide homological information.

\section{Preliminaries}

\label{section:preliminaryresults}We begin with a list of alternative
characterizations of the class of real number fields $F$ of interest. To fix
terminology, recall that the $r$ real embeddings of $F$ define the
\emph{signature} map $F^{\times}/(F^{\times})^{2}\rightarrow(\mathbb{Z}%
/2)^{r}$ where a unit is mapped to the signs of its images under the real
embeddings. One says that $R_{F}=\mathcal{O}_{F}[\frac{1}{2}]$, the ring of
$2$-integers in $F$, has \emph{units of independent signs} if the signature
map remains surjective when restricted to the square classes $R_{F}^{\times
}/(R_{F}^{\times})^{2}$ of $R_{F}^{\times}$. A \emph{dyadic prime} in $F$ is a
prime ideal in the ring of integers $\mathcal{O}_{F}$ lying over the rational
prime ideal $(2)$. The \emph{narrow Picard group} \textrm{Pic}$_{+}(R_{F})$
consists of fractional $R_{F}$-ideals modulo totally positive principal
ideals, defined as in \cite[V\S 1]{FrohlichTaylor}.

The \emph{Witt ring} $W(A)$ of a commutative unital ring $A$ with involution
is defined in terms of Witt classes of symmetric nondegenerate bilinear forms
\cite[I (7.1)]{MH}; in $K$-theoretic terms, as an abelian group it coincides
with the cokernel ${}_{1}W_{0}(A)$ of the hyperbolic map $K_{0}(A)\rightarrow
{}_{1}KQ_{0}(A)$ if $2$ is invertible in $A$ \cite{K:AnnM112fun}.
Symmetrically, we define the \emph{coWitt group} $W^{\prime}(A)={}_{1}%
W_{0}^{\prime}(A)$ as the kernel of the forgetful rank map ${}_{1}%
KQ_{0}(A)\rightarrow{}K_{0}(A)$. In this section especially, we often use the
simplified notations $W(A)$ and $W^{\prime}(A)$ instead of ${}_{1}W_{0}(A)$
and $_{1}W_{0}^{\prime}(A)$.

In the statement below, for a finite abelian group $G$, we write $G\{2\}$ for
its $2$-Sylow subgroup.

\begin{proposition}
\label{2+-regular characterization} Let $F$ be a number field with $r$ real
embeddings and $c$ pairs of complex embeddings. Then the following are equivalent.

\begin{enumerate}
\item $F$ is $2$-regular; that is, the $2$-Sylow subgroup of the finite
abelian group $K_{2}(\mathcal{O}_{F})$ has order $2^{r}$.

\item The real embeddings of $F$ induce isomorphisms on $2$-Sylow subgroups
\[
K_{2}(\mathcal{O}_{F})\{2\}\overset{\cong}{\longrightarrow}K_{2}%
(R_{F})\{2\}\overset{\cong}{\longrightarrow}\oplus^{r}K_{2}(\mathbb{R}%
)\cong(\mathbb{Z}/2)^{r}.
\]

\item The finite abelian group $K_{2}(\mathcal{O}_{F(\sqrt{-1})})$ has odd order.

\item There is a unique dyadic prime in $F$ and the narrow Picard group
\textrm{Pic}$_{+}(R_{F})$ has odd order.

\item There is a unique dyadic prime in $F$, the Picard group \textrm{Pic}%
$(R_{F})$ has odd order, and $R_{F}$ has units with independent signs.

\item The nilradical of the Witt ring $W(R_{F})$ is a finite abelian group of
order $2^{c+1}$. In particular, if $F$ is a totally real number field, so that
$c=0$, then the nilradical of $W(R_{F})$ has order $2$.

\item The Witt ring $W(R_{F})$ is a finitely generated abelian group of rank
$r$ with torsion subgroup of order $2^{c+1}$. In particular, if $F$ is a
totally real number field, then $W(R_{F})$ is isomorphic to $\mathbb{Z}%
^{r}\oplus\mathbb{Z}/2$.
\end{enumerate}

If $F$ is a totally real number field and satisfies any of the equivalent
conditions above, then the free part of $W(R_{F})$ is generated by elements
$\left\langle 1\right\rangle ,\left\langle u_{1}\right\rangle ,\dots
,\left\langle u_{r-1}\right\rangle $ where $u_{i}\in R_{F}^{\times}$ is
negative at the $i$\thinspace th embedding and positive elsewhere.
\end{proposition}

\noindent\textbf{Proof.} Obviously, (2) implies (1). For the converse, we
apply Tate's $2$-rank formula for $K_{2}$ \cite[Theorem 6.2]{Tate}, which
shows that the group $K_{2}(\mathcal{O}_{F})$ maps, via square power norm
residue symbols, onto the subgroup of the direct sum (over all archimedean
places and places over $2$) of copies of $\mu_{2}(F)=\mathbb{Z}/2$ that
consists of the elements $z=(z_{v})$ such that $\sum z_{v}=0$. This forces $r$
to be a lower bound for the $2$-rank of $K_{2}(\mathcal{O}_{F})$. Of course,
the unique group of order $2^{r}$ and $2$-rank at least $r$ is $(\mathbb{Z}%
/2)^{r}$; so, the converse follows.

See \cite[Proposition 2.2]{RO} for the equivalence between (2) and (4). By
\cite[(4.1),(4.6)]{ConnerHurrelbrink}, again using \cite[Theorem 6.2]{Tate},
(2), (3) and (5) are equivalent. The equivalences between (4), (6), and (7)
are immediate from \cite[Corollary 3.6, Theorem 4.7]{Czogala}. Finally, the
claim concerning the generators of the free part of $W(R_{F})$ follows as in
\cite[IV (4.3)]{MH}. \hfill$\Box$

\begin{remark}
1. Further equivalent conditions, in terms of \'{e}tale cohomology, appear in
\cite[Proposition 2.2]{RO}.

2. Berger \cite{Berger} calls a totally real number field $F$ satisfying (5)
above $2^{+}$\emph{-regular}, and has shown in \cite{Berger} that each totally
real $2$-regular number field has infinitely many totally real quadratic field
extensions satisfying the equivalent conditions in Proposition
\ref{2+-regular characterization}. In particular, there exist totally real
$2$-regular number fields of arbitrarily high degree.

3. In the other direction, by \cite[(4.1)]{ConnerHurrelbrink} every subfield
of a totally real $2$-regular number field also enjoys this property.

4. In \cite[pg.~95]{MH}, it is shown that for totally real $F$ there is also
an equivalence between amended conditions (4)--(7) above, in which $R_{F}$ is
replaced by $\mathcal{O}_{F}$ and $2^{c+1}$ by $2^{c}$. However, in view of
the surjection \textrm{Pic}$_{+}(\mathcal{O}_{F})\twoheadrightarrow
$\textrm{Pic}$_{+}(R_{F})$, these amended conditions define a strict subclass
of those considered here. After Gauss, one knows that for a real quadratic
number field $F=\mathbb{Q}(\sqrt{d})$ the narrow Picard group \textrm{Pic}%
$_{+}(\mathcal{O}_{F})$ is an odd torsion group if and only if $F$ has
prime-power discriminant; that is, $d=2$ or $d=p$ with $p\equiv
1\;(\mathrm{mod}\ 4)$ a prime number. Thus $\mathbb{Q}(\sqrt{p})$ with
$p\equiv3\;(\mathrm{mod}\ 4)$ a prime number and $\mathbb{Q}(\sqrt{2p})$ with
$p\equiv\pm3\;(\mathrm{mod}\ 8)$ a prime number fail to lie in the subclass,
cf.~\cite{CO}.
\end{remark}

\textit{Suppose henceforth} that $2$ is a unit in a domain $A$ of dimension at
most $1$, \textsl{e.g}.~a subring of some number field $F$. Let $k_{0}(A)$
denote the $0$\thinspace th Tate cohomology group $\widehat{H}^{0}%
(\mathbb{Z}/2;\,K_{0}(A))$ of $\mathbb{Z}/2$ acting on $K_{0}(A)$ by duality
as in \cite[pg.~278]{K:AnnM112fun}, and $k_{0}^{\prime}(A)$ denote the
$1$\thinspace st Tate cohomology group $\widehat{H}^{1}(\mathbb{Z}%
/2;\,K_{0}(A))$. \vspace{0.1in}

There is a well defined induced \emph{rank map }$\rho\colon W(A)\rightarrow
k_{0}(A)$, whose image always has a $\mathbb{Z}/2$ summand. From the 12-term
exact sequence established in \cite[pg.~278]{K:AnnM112fun} there is an exact
sequence
\[
k_{0}^{\prime}(A)\longrightarrow W^{\prime}(A)\overset{\varphi}%
{\longrightarrow}W(A)\overset{\rho}{\longrightarrow}k_{0}(A),
\]
where $\varphi$ is simply the composition $W^{\prime}(A)\hookrightarrow{}%
_{1}KQ_{0}(A)\twoheadrightarrow W(A)$.

The next result is almost immediate, but worth recording. Part (c) uses the
observation that $\tilde{K}_{0}(A)\cong\mathrm{Pic}(A)$ when $A$ is of
dimension at most $1$. If further \textrm{Pic}$(A)$ is an odd torsion group,
then $k_{0}^{\prime}(A)=0$ and $k_{0}(A)\cong\mathbb{Z}/2$. Also, (d) refers
to the fact that, when $A$ is a field, its \emph{fundamental ideal} is the
unique ideal $I$ of $W(A)$ with $W(A)/I\cong\mathbb{F}_{2}$ \cite[III
(3.3)]{MH}.

\begin{lemma}
\label{proposition:dedekindandoddtorsionpicardgroup}(a) In general, $\varphi$
maps the coWitt group of $A$ onto the kernel of the rank map $\rho$.

(b) When $k_{0}^{\prime}(A)=0$, then $\varphi$ is an isomorphism between
$W^{\prime}(A)$ and $\mathrm{Ker}\rho$.

(c) When $A$ has dimension at most $1$ and $\mathrm{Pic}(A)$ is an odd torsion
group, then there is a natural short exact sequence%
\[
0\rightarrow W^{\prime}(A)\overset{\varphi}{\longrightarrow}W(A)\overset{\rho
}{\longrightarrow}\mathbb{Z}/2\rightarrow0\text{.}%
\]

(d) When $A$ is a field, then $\varphi$ induces an isomorphism between the
coWitt group of $A$ and the fundamental ideal of the Witt ring.$\hfill\Box$
\end{lemma}

\vspace{0.1in}

We can now give a further characterization of the class of $2$-regular totally
real number fields, in terms of the coWitt group.

\begin{lemma}
\label{2-regular coWitt}A totally real number field $F$ with $r$ real
embeddings is $2$-regular if and only if both

\begin{enumerate}
\item[(i)] the Picard group $\mathrm{Pic}(R_{F})$ has odd order, and

\item[(ii)] the coWitt group $W^{\prime}(R_{F})$ is isomorphic to
$\mathbb{Z}^{r}\oplus\mathbb{Z}/2$.
\end{enumerate}
\end{lemma}

\noindent\textbf{Proof.} Here we use the fact that, as for both $K_{0}$ and
$KQ_{0}$, the Witt group $W(\mathbb{R})=W_{0}(\mathbb{R})$ and the coWitt
group $W^{\prime}(\mathbb{R})=W_{0}^{\prime}(\mathbb{R})$ do not depend on the
topology of $\mathbb{R}$.

In both directions of the statement of the lemma (one way uses Proposition
\ref{2+-regular characterization}(5)), we have from Lemma
\ref{proposition:dedekindandoddtorsionpicardgroup}(c) the map of short exact
sequences (where the middle vertical map is surjective):
\[%
\begin{array}
[c]{ccccccc}%
0\rightarrow & W^{\prime}(R_{F}) & \overset{\varphi}{\longrightarrow} &
W(R_{F}) & \overset{\rho}{\longrightarrow} & \mathbb{Z}/2 & \rightarrow0\\
& \downarrow &  & \downarrow &  & \downarrow\mathrm{id} & \\
0\rightarrow & W^{\prime}(\mathbb{R}) & \overset{\varphi}{\longrightarrow} &
W(\mathbb{R}) & \overset{\rho}{\longrightarrow} & \mathbb{Z}/2 & \rightarrow0
\end{array}
\]
Since $W(\mathbb{R})\cong\mathbb{Z}$ \cite[III (2.7)]{MH}, in the upper
sequence $\rho$ must be trivial on torsion elements. By combining with
Proposition \ref{2+-regular characterization}(7), we obtain the result.
\hfill$\Box$

When $F$ has $r$ real embeddings, for $A=R_{F},F$, the map $W(A)\rightarrow
W(\mathbb{R})^{r}\cong\mathbb{Z}^{r}$ is the \emph{total signature }$\sigma$
\cite[III (2.9)]{MH}. \vspace{0.1in}

To define further invariants, we recall from \cite{K:AnnM112fun}
generalizations of some of the definitions above. For $\varepsilon=\pm1$, and
$n\geq1$, we set
\begin{align*}
{}_{\epsilon}W_{n}(A)  &  =\mathrm{Coker}[K_{n}(A)\rightarrow{}_{\epsilon
}KQ_{n}(A)],\\
{}_{\epsilon}W_{n}^{\prime}(A)  &  =\mathrm{Ker}[{}_{\epsilon}KQ_{n}%
(A)\rightarrow K_{n}(A)],\\
k_{n}(A)  &  =\widehat{H}^{0}(\mathbb{Z}/2;\,K_{n}(A)),\\
k_{n}^{\prime}(A)  &  =\widehat{H}^{1}(\mathbb{Z}/2;\,K_{n}(A)).
\end{align*}

The next invariant we shall employ is the \emph{discriminant map }%
${}_{\epsilon}W_{0}^{\prime}(A)\rightarrow k_{1}^{\prime}(A)$. To recall the
definition, suppose that $M$ and $N$ are quadratic modules with isomorphic
underlying $A$-modules. The elements of ${}_{\epsilon}W_{0}^{\prime}(A)$ are
of the form $M-N$ with $M$ and $N$ isomorphic modules. An isomorphism
$\alpha\colon M\rightarrow N$ induces an automorphism $\alpha^{\ast}\alpha$ of
$M$ that is antisymmetric for the $\mathbb{Z}/2$-action. Its class in
$k_{1}^{\prime}(A)$, which is independent of $\alpha$, defines the desired
invariant. For a generalization of the above we refer to \cite{K:AnnM112fun}.

If $SK_{1}(A)=0$ (\textsl{e.g. }$A$ a ring of $S$-integers in a number field
\cite[Corollary 4.3]{BassMilnorSerre} or $A$ a field), then $k_{1}^{\prime
}(A)$ is isomorphic to the group of square classes $A^{\times}/(A^{\times
})^{2}$ of units in $A$ and $k_{1}(A)=\{\pm1\}$. This affords the following description.

\begin{lemma}
\label{SK_1(A) = 0}Assume that $SK_{1}(A)=0$\textsf{.}

Then the rank/determinant map ${}_{1}W_{1}(A)\rightarrow k_{1}(A)=\{\pm1\}$ is
surjective, and there is a short exact sequence%
\[
0\rightarrow{}_{-1}W_{2}(A)\longrightarrow{}_{1}W_{0}^{\prime}(A)\overset
{\mathrm{disc}}{\longrightarrow}k_{1}^{\prime}(A)=A^{\times}/(A^{\times}%
)^{2}\rightarrow0\text{.}%
\]

If $A$ is a field, we may identify $W^{\prime}(A)={}_{1}W_{0}^{\prime}(A)$
(resp.~ ${}_{-1}W_{2}(A)$) with the fundamental ideal in the Witt group $W(A)$
(resp.~ its square), so that the exact sequence above reduces to
\[
0\rightarrow I^{2}\longrightarrow I\longrightarrow I/I^{2}\rightarrow0\text{.}%
\]
The result is still true for the field $\mathbb{R}$ and $\mathbb{C}$ with
their usual topology.
\end{lemma}

\noindent\textbf{Proof.} Since $SK_{1}(A)=0$, we have $K_{1}(A)=A^{\times}$.
Therefore, $k_{1}(A)$ is reduced to $\pm1$ which is a $1\times1$ unitary
matrix. This explains the surjectivity. By the exact sequence
\[
{}_{1}W_{1}(A)\longrightarrow k_{1}(A)\longrightarrow{}_{-1}W_{2}%
(A)\longrightarrow{}_{1}W_{0}^{\prime}(A)\longrightarrow k_{1}^{\prime}(A)
\]
from \cite{K:AnnM112fun}, ${}_{-1}W_{2}(A)$ now identifies with the kernel of
the induced discriminant map $W^{\prime}(A)\rightarrow k_{1}^{\prime}(A)\cong
A^{\times}/(A^{\times})^{2}$.

Moreover, taking $\varepsilon=1$, $W^{\prime}(A)$ surjects onto $A^{\times
}/(A^{\times})^{2}$ since the discriminant maps $\left\langle u\right\rangle
-\left\langle 1\right\rangle $ to $u\in A^{\times}/(A^{\times})^{2}$.

When $A$ is a field, by Lemma
\ref{proposition:dedekindandoddtorsionpicardgroup}\thinspace(d) above the
coWitt group identifies with the fundamental ideal $I$, while by \cite[III
(5.2)]{MH}\textsf{ }the discriminant induces an isomorphism $I/I^{2}\cong
A^{\times}/(A^{\times})^{2}$. Finally, if $A=\mathbb{R}$ or $\mathbb{C}$ with
the usual topology, it is easy to see that the group${}_{-1}W_{2}(A)$
coincides with the same group when we regard $\mathbb{R}$ and $\mathbb{C}$
with the discrete topology, thanks to the fundamental theorem in topological
hermitian $K$-theory proved in \cite{Karoubi LNM343}.\hfill$\Box$\smallskip

\begin{examples}
\label{examples for A/A^{2}}In the case of $A=\mathbb{R}$, from \cite[III
(2.7)]{MH} the discriminant map above is the surjection $\mathbb{Z}%
\rightarrow\mathbb{Z}/2$.

For $A=\mathbb{F}_{q}$ by \cite[III (5.2), (5.9)]{MH} it is the isomorphism
$A^{\times}/(A^{\times})^{2}\cong\mathbb{Z}/2$.

For $A=R_{F}$, with $F$ totally real, the Dirichlet $S$-unit theorem for
$R_{F}$ \cite[Ch. IV Theorem 9]{Weil}\textsf{ }gives $R_{F}^{\times}%
\cong\mathbb{Z}^{r+d-1}\times\mu(F)$, where $d$ is the number of dyadic
places, and the group of roots of unity $\mu(F)$ has order $2$. In the
$2$-regular case, $d=1$ by Proposition \ref{2+-regular characterization}(4),
so $R_{F}^{\times}/(R_{F}^{\times})^{2}\cong(\mathbb{Z}/2)^{r+1}$.
\end{examples}

What follows is a key ingredient in the proof of Theorem \ref{theorem1}.

\begin{proposition}
\label{proposition:keyinput} Suppose that $F$ is a totally real $2$-regular
number field with $r$ real embeddings. Then the residue field map
$R_{F}\rightarrow\mathbb{F}_{q}$ and the real embeddings of $F$ induce an
isomorphism between coWitt groups
\[%
\begin{array}
[c]{c}%
W^{\prime}(R_{F})\overset{\cong}{\longrightarrow}\bigoplus\limits^{r}%
W^{\prime}(\mathbb{R})\oplus W^{\prime}(\mathbb{F}_{q})\cong\mathbb{Z}%
^{r}\oplus\mathbb{Z}/2\text{.}%
\end{array}
\]

\end{proposition}

\noindent\textbf{Proof.} Recalling that $SK_{1}$ is trivial for $R_{F}$ and
$F$ \cite[Corollary 4.3]{BassMilnorSerre}, we obtain from the last lemma a map
of short exact sequences
\[%
\begin{array}
[c]{ccccccc}%
0\rightarrow & {}_{-1}W_{2}(R_{F}) & \rightarrow & W^{\prime}(R_{F}) &
\rightarrow & R_{F}^{\times}/(R_{F}^{\times})^{2} & \rightarrow0\\
& \downarrow &  & \downarrow &  & \downarrow & \\
0\rightarrow & \bigoplus\limits^{r}{}_{-1}W_{2}(\mathbb{R})\oplus{}_{-1}%
W_{2}(\mathbb{F}_{q}) & \rightarrow & \bigoplus\limits^{r}W^{\prime
}(\mathbb{R})\oplus W^{\prime}(\mathbb{F}_{q}) & \rightarrow & \bigoplus
\limits^{r}\mathbb{R}^{\times}/(\mathbb{R}^{\times})^{2}\oplus\mathbb{F}%
_{q}^{\times}/(\mathbb{F}_{q}^{\times})^{2} & \rightarrow0
\end{array}
\]
which by Lemma \ref{2-regular coWitt} and Example (\ref{examples for A/A^{2}})
above\textsf{ }takes the form:
\begin{equation}%
\begin{array}
[c]{ccccccc}%
0\rightarrow & \mathbb{Z}^{r} & \rightarrow & \mathbb{Z}^{r}\oplus
\mathbb{Z}/2 & \rightarrow & (\mathbb{Z}/2)^{r+1} & \rightarrow0\\
& \downarrow &  & \downarrow &  & \downarrow & \\
0\rightarrow & \bigoplus\limits^{r}\mathbb{Z}\oplus0 & \rightarrow &
\bigoplus\limits^{r}\mathbb{Z}\oplus\mathbb{Z}/2 & \rightarrow &
\bigoplus\limits^{r}\mathbb{Z}/2\oplus\mathbb{Z}/2 & \rightarrow0
\end{array}
\label{diagram: computed discriminant}%
\end{equation}

From the homotopy cartesian square (\ref{algebraic K Bokstedt square}) of
$2$-completed spectra, established in \cite{HO} and \cite{Mitchell} -- see
Appendix A, we deduce the short exact sequence of $2$-completed groups
\begin{equation}%
\begin{array}
[c]{c}%
0\rightarrow\bigoplus^{r}K_{2}(\mathbb{C})\rightarrow K_{1}(R_{F}%
)\rightarrow\bigoplus^{r}K_{1}(\mathbb{R})\oplus K_{1}(\mathbb{F}%
_{q})\rightarrow0\text{.}%
\end{array}
\end{equation}
(The left hand $0$ is justified by the facts that $K_{2}(\mathbb{F}_{q})=0$
and $K_{2}(\mathbb{R})=\mathbb{Z}/2$, while $K_{2}(\mathbb{C})$ is
torsion-free\textsf{.) }Therefore, the final vertical map in
(\ref{diagram: computed discriminant}) is an isomorphism. Observe that
commutativity of the right-hand square implies that the middle vertical map in
(\ref{diagram: computed discriminant}) is a monomorphism on the\textsf{
}$\mathbb{Z}/2$ summand.\textsf{ }Thus, to show that the middle vertical map
is an isomorphism, because this is equivalent to its being an
epimorphism,\textsf{ }it suffices to show that either of the homomorphisms
$W^{\prime}(R_{F})\rightarrow\bigoplus\limits^{r}W^{\prime}(\mathbb{R})$ or
${}_{-1}W_{2}(R_{F})\rightarrow\bigoplus\limits^{r}{}_{-1}W_{2}(\mathbb{R})$
is surjective.

For the former homomorphism, we use the map of short exact sequences afforded
by Lemma \ref{proposition:dedekindandoddtorsionpicardgroup}(c):%
\[%
\begin{array}
[c]{ccccccc}%
0\rightarrow & W^{\prime}(R_{F}) & \longrightarrow & W(R_{F}) &
\longrightarrow & \mathbb{Z}/2 & \rightarrow0\\
& \downarrow\sigma^{\prime} &  & \downarrow\sigma &  & \downarrow\Delta & \\
0\rightarrow & \bigoplus\limits^{r}W^{\prime}(\mathbb{R}) & \longrightarrow &
\bigoplus\limits^{r}W(\mathbb{R}) & \longrightarrow & \bigoplus\limits^{r}%
\mathbb{Z}/2 & \rightarrow0
\end{array}
\]
Since the diagonal map $\Delta$ is injective, the snake lemma gives the exact
sequence%
\[
0\longrightarrow\mathrm{Coker}\sigma^{\prime}\longrightarrow\mathrm{Coker}%
\sigma\longrightarrow(\mathbb{Z}/2)^{r-1}\longrightarrow0\text{.}%
\]
According to \cite[Corollary 4.8]{Czogala} (applicable because by Proposition
\ref{2+-regular characterization}(4) \textrm{Pic}$_{+}(R_{F})$ has odd order),
$\sigma:W(R_{F})\rightarrow\bigoplus\limits^{r}W(\mathbb{R})$ has the same
cokernel as $\sigma:W(F)\rightarrow\bigoplus\limits^{r}W(\mathbb{R})$. Now, by
\cite[pp.~64-65]{MH} (applicable because by Proposition
\ref{2+-regular characterization}(5) $R_{F}$ has units with independent
signs), $\mathrm{Coker}\sigma$ is $(\mathbb{Z}/2)^{r-1}$. Hence,
$\mathrm{Coker}\sigma^{\prime}=0$, as desired.

The second approach, showing that ${}_{-1}W_{2}(R_{F})\rightarrow
\bigoplus\limits^{r}{}_{-1}W_{2}(\mathbb{R})$ is surjective, instead of
\cite{Czogala} uses the \textquotedblleft classical\textquotedblright%
\ treatment of \cite{MH}, by invoking the generalized Hasse-Witt invariant for
quadratic forms.

Recall the map ${}_{-1}W_{2}(A)\longrightarrow k_{2}(A)$ employed in the
definition of the 12-term exact sequence in \cite{K:AnnM112fun}. If $A$ is a
field, then by Lemma \ref{SK_1(A) = 0} above we may identify ${}_{-1}W_{2}(A)$
with $I^{2}(A)$; according to \cite{Guin}, this map from $I^{2}(A)$ to
$k_{2}(A)=K_{2}(A)/2$ gives an equivalent definition of the classical
\emph{Hasse-Witt invariant}. Moreover, the Hasse-Witt invariants for $R_{F}$
and $F$ induce a commutative diagram with exact rows:
\[%
\begin{array}
[c]{ccccc}%
{}_{-1}W_{2}(R_{F}) & \rightarrow & k_{2}(R_{F}) & \rightarrow & 0\\
\downarrow &  & \downarrow &  & \\
{}_{-1}W_{2}(F) & \rightarrow & k_{2}(F) & \rightarrow & 0
\end{array}
\]
Now, from \cite{HO} and \cite{Mitchell}, the real embeddings of $F$ induce a
surjective map from $R_{F}^{\times}$ to $\oplus^{r}K_{1}(\mathbb{R})$. Hence,
${}_{-1}W_{2}(R_{F})$ surjects onto $k_{2}(R_{F})\cong(\mathbb{Z}/2)^{r}$
because the tensor product $(\left\langle u\right\rangle -\left\langle
1\right\rangle )\otimes(\left\langle v\right\rangle -\left\langle
1\right\rangle )$ maps to the symbol $\{u,v\}$ in $k_{2}(R_{F})$. By reference
to \cite[III (5.9)]{MH}, where the arguments can be extended to any ring $D$
of $S$-integers in a number field, the kernels of the Hasse-Witt surjections
for $R_{F}$ and $F$ are isomorphic to $I^{3}(F)\cong\bigoplus^{r}%
I^{3}(\mathbb{R})\cong8\mathbb{Z}^{r}$ via the signature map, and hence there
is a map of short exact sequences:%
\[%
\begin{array}
[c]{ccccccc}%
0\rightarrow & 8\mathbb{Z}^{r} & \longrightarrow & {}_{-1}W_{2}(R_{F}) &
\longrightarrow & (\mathbb{Z}/2)^{r} & \rightarrow0\\
& \downarrow\cong &  & \downarrow &  & \downarrow\cong & \\
0\rightarrow & 8\mathbb{Z}^{r} & \longrightarrow & \oplus^{r}{}_{-1}%
W_{2}(\mathbb{R}) & \overset{\chi}{\longrightarrow} & (\mathbb{Z}/2)^{r} &
\rightarrow0
\end{array}
\]
where the notation $8\mathbb{Z}^{r}$ means the image in $\oplus^{r}%
W(\mathbb{R})=\mathbb{Z}^{r}$ of the kernel of $\chi$ by the total signature
map. The $5$-lemma now shows that ${}_{-1}W_{2}(R_{F})$ is isomorphic to
$\oplus^{r}{}_{-1}W_{2}(\mathbb{R})$ by the middle vertical map. \hfill$\Box$

\bigskip

We finish this section by relating the numbers $w_{m}$ in the formulation of
Theorem \ref{theorem2} to $t_{n}$, the $2$-adic valuation $(q^{(n+1)/2}%
-1)_{2}$ of $q^{(n+1)/2}-1$ for $n$ odd.

\begin{lemma}
Suppose that $q\equiv1\;(\mathrm{mod}\ 4)$, and write $(-)_{2}$ for the
$2$-adic valuation. Then $(q^{m}-1)_{2}=(q-1)_{2}(m)_{2}$.
\end{lemma}

\noindent\textbf{Proof.} With $q=4r+1$ and $t=(m)_{2}$, note that $q^{t}-1$ is
divisible by $4rt$ but not by $8rt$, due to the binomial identity. Set
$s=8t(r)_{2}$ and $u=m/t$. Since $(\mathbb{Z}/2^{s})^{\times}$ has even order
and $u$ is odd, $q^{m}-1=(q^{t})^{u}-1$ is divisible by $s/2$ but not by $s$.
\hfill$\Box$

\begin{lemma}
Let $q$ be an odd number. If $m$ is odd, then $(q^{m}-1)_{2}=(q-1)_{2}$. If
$m=2m^{\prime}$ is even, then $(q^{m}-1)_{2}=(q^{2}-1)_{2}(m^{\prime})_{2}$.
\end{lemma}

\noindent\textbf{Proof.} If $m$ is odd, writing $q^{m}-1=(q-1)(q^{m-1}%
+\cdots+1)$ shows that $(q^{m}-1)_{2}=(q-1)_{2}$ since $(q^{m-1}+\cdots+1)$ is
odd. If $m=2m^{\prime}$ is even, write $q^{m}-1$ as $(q^{2})^{m^{\prime}}-1$
where $q^{2}\equiv1\;(\mathrm{mod}\ 4)$, and apply the previous lemma.
\hfill$\Box$

\medskip

Returning now to the setting of this paper, we have a prime $q$ defined, as in
the Introduction, in terms of the number field $F$ such that $q$ is $\equiv
\pm1\;(\mathrm{mod}\ 2^{a_{F}})$ but not $(\mathrm{mod}\ 2^{a_{F}+1})$. Thus,
for even $m$, the last lemma above gives
\[
(q^{m}-1)_{2}=2^{a_{F}}(m)_{2}=:w_{m}\text{,}%
\]
where the number $w_{m}$ appears in Theorem \ref{theorem2}. Hence, we have the
following description of $t_{n}:=(q^{(n+1)/2}-1)_{2}$ in terms of $w_{m}$.

\begin{lemma}
\label{lemma: t and w}Let $F$ be a totally real number field for which $q$ is
chosen as in the Introduction. Then, for $n\equiv3\;(\mathrm{mod}\ 4)$,%
\[
t_{n}=w_{(n+1)/2}\text{.}\vspace{-20pt}%
\]
\hfill$\Box\smallskip$
\end{lemma}

\section{Proof of Theorem \ref{theorem1}}

\label{section:proofoftheorem1} Before commencing consideration of Theorem
\ref{theorem1}, we make a few preparatory remarks.

The following result, of independent interest, will be used below.

\begin{proposition}
\label{Wodd}Let $R$ be any ring of $S$-integers in a number field $F$ with
$1/2\in R$. Then the odd torsion of the generalized Witt groups%
\[
{}_{\varepsilon}W_{n}(R)=\mathrm{Coker}\left[  K_{n}(R)\longrightarrow
{}_{\varepsilon}KQ_{n}(R)\right]
\]
and the coWitt groups%
\[
{}_{\varepsilon}W_{n}^{^{\prime}}(R)=\mathrm{Ker}\left[  {}_{\varepsilon
}KQ_{n}(R)\longrightarrow K_{n}(R)\right]
\]
is trivial for all values of $n\in%
\mathbb{Z}%
$.
\end{proposition}

\noindent\textbf{Proof.} We first remark that the higher Witt groups and
coWitt groups have the same odd torsion because of the $12$ term exact
sequence proved in \cite{K:AnnM112fun}. The same exact sequence shows that we
have an isomorphism modulo odd torsion between ${}_{\varepsilon}W_{n}(R)$ and
${}_{-\varepsilon}W_{n+2}(R)$. Therefore, we have only to compute the four
cases $\varepsilon=\pm1$ and $n=0,1$.

For $\varepsilon=-1$, recall that ${}_{-1}KQ_{1}(R)=\mathrm{Sp}(R)\left/
\left[  \mathrm{Sp}(R),\,\mathrm{Sp}(R)\right]  \right.  =0$ according to
\cite[Corollary $4.3$]{BassMilnorSerre}, and therefore its quotient ${}%
_{-1}W_{1}(R)=0$. By \cite[I (3.5)]{MH} there are isomorphisms ${}_{-1}%
KQ_{0}(R)\cong\mathbb{Z}$ detected by the (even) rank of the free symplectic
$R$-inner product space. Therefore, ${}_{-1}W_{0}(R)=\mathrm{Coker}\left[
K_{0}(R)\longrightarrow{}_{-1}KQ_{0}(R)\right]  =0.$

For $\varepsilon=1,$ it is well known (see for instance \cite[Corollary
IV.3.3]{MH}) that the Witt group ${}_{1}W(R)$ injects in ${}_{1}W(F)$ and that
the torsion of ${}_{1}W(F)$ is $2$-primary according to the same book
\cite[Theorem III.3.10]{MH}. Finally, we consider the exact sequence%
\[
K_{1}(R)\longrightarrow{}_{1}KQ_{1}(R)\longrightarrow{}_{2}\mathrm{Pic}%
(R)\oplus Z_{2}(R)\longrightarrow0
\]
proved in \cite[(4.7.6)]{Bass} with different notations. Here, ${}%
_{2}\mathrm{Pic}(R)$ is the $2$-torsion of the Picard group and $Z_{2}(R)$ is
the group of locally constant maps from $\mathrm{Spec}(R)$ to $%
\mathbb{Z}
/2$. Moreover, $\mathrm{Spec}(R)$ is connected, so that $Z_{2}(R)=%
\mathbb{Z}
/2$, and the previous exact sequence yields the isomorphism%
\begin{equation}
{}_{1}W(R)\cong{}_{2}\mathrm{Pic}(R)\oplus%
\mathbb{Z}
/2\text{.}\vspace{-10pt}\label{W1}%
\end{equation}
\hfill$\Box\smallskip$

The next follows by comparing the splittings of the hermitian $K$-theory
spectrum and the $K$-theory spectrum according to their canonical involutions
\cite[pg.~253]{K:AnnM112fun}. We record it for the sake of completeness: by
the previous result, it may be applied to any ring of $2$-integers in a number field.

\begin{proposition}
Let $A$ be any ring, with $n$ an integer such that the odd torsion
${}_{\varepsilon}W_{n}(A)$ vanishes. Then the odd torsion subgroup of
${}_{\varepsilon}KQ_{n}(A)$ is the invariant part of the odd torsion subgroup
of $K_{n}(A)$ induced by the involution $M\mapsto{}^{t}M^{-1}$ on
$\mathrm{GL}(A)$.\hfill$\Box\smallskip$
\end{proposition}

We now turn to consideration of spectra. The following simple exercise in
stable homotopy theory ((iv), (v) are from \cite[Chapter VI (5.1),(5.2)]%
{BK304}, suitably adopted to spectra in \cite[Proposition 2.5]{Bousfieldloc})
will be used implicitly at various times. Here and subsequently, we write
$\Omega\mathcal{E}$ to indicate the shifted spectrum whose space at level $n$
is $\mathcal{E}_{n-1}\simeq\Omega(\mathcal{E}_{n})$; similarly for other
shifts. In (iv), $\mathbb{Z}_{2^{\infty}}$ denotes the quasicyclic $2$-group
$\underrightarrow{\mathrm{lim}}\,\mathbb{Z}/2^{n}$.

\begin{lemma}
\label{Connective.Spectra} For any spectrum $\mathcal{E}$ and fiber sequence
$\mathcal{F}\rightarrow\mathcal{E}\rightarrow\mathcal{B}$ of spectra, the
connective covers and $2$-completions have the following properties.

\begin{enumerate}
\item[\emph{(i)}] $\mathcal{F}^{c}\rightarrow\mathcal{E}^{c}\rightarrow
\mathcal{B}^{c}$ is a fiber sequence, provided that $\pi_{0}(\mathcal{E}%
)\rightarrow\pi_{0}(\mathcal{B})$ is an epimorphism.

\item[\emph{(ii)}] $(\Omega\mathcal{E})^{c}\simeq\Omega(\mathcal{E}^{c})$,
provided that $\pi_{0}(\mathcal{E})=0$.

\item[\emph{(iii)}] $(\Omega^{-1}\mathcal{E})^{c}\simeq\Omega^{-1}%
(\mathcal{E}^{c})$, provided that $\pi_{-1}(\mathcal{E})=0$.

\item[\emph{(iv)}] $(\mathcal{E}_{\#})^{c}\simeq(\mathcal{E}^{c})_{\#}%
$\thinspace, provided that $\mathrm{Hom}(\mathbb{Z}_{2^{\infty}},$%
\thinspace$\pi_{-1}(\mathcal{E}))=0$, as for example when $\pi_{-1}%
(\mathcal{E})$ is finitely generated.

\item[\emph{(v)}] When all homotopy groups of $\mathcal{E}$ are finitely
generated, then for all $i\in\mathbb{Z}$ the $2$-completed spectrum
$\mathcal{E}_{\#}$ has
\[
\pi_{i}(\mathcal{E}_{\#})=\pi_{i}(\mathcal{E})\otimes\mathbb{Z}_{2}%
\text{.}\vspace{-20pt}%
\]
\hfill$\Box\smallskip$
\end{enumerate}
\end{lemma}

Even though (iv) above applies to our spectra whose homotopy groups are
finitely generated, we clarify our convention by defining $\mathcal{E}%
_{\#}^{c}$ as $(\mathcal{E}^{c})_{\#}$. This spectrum $\mathcal{E}_{\#}^{c}$
is named the $2$-completed connective spectrum associated to $\mathcal{E}$.
For instance,%
\[
\mathcal{K}(R_{F})_{\#}^{c}=(\mathcal{K}(R_{F})^{c})_{\#}=\mathcal{K}%
(R_{F})_{\#}\,\text{,}%
\]
and for $n\geq0$%
\[
\pi_{n}(\mathcal{K}(R_{F})_{\#}^{c})=K_{n}(R_{F})\otimes\mathbb{Z}%
_{2}\,\text{;}%
\]
while
\[
{}_{\varepsilon}\mathcal{KQ}(R_{F})_{\#}^{c}=({}_{\varepsilon}\mathcal{KQ}%
(R_{F})^{c})_{\#}\,\text{,}%
\]
and for $n\geq0$%
\[
\pi_{n}({}_{\varepsilon}\mathcal{KQ}(R_{F})_{\#}^{c})={}_{\varepsilon}%
KQ_{n}(R_{F})\otimes\mathbb{Z}_{2}\,\text{.}%
\]
This uses the result of \cite[(3.6), pg. 795.]{BK} that the groups
${}_{\varepsilon}KQ_{n}(R_{F})$ are finitely generated. For completeness'
sake, we give an alternative proof in Proposition \ref{Finiteness} at the end
of this section.

\medskip The following lemma is used later on in our argument.

\begin{lemma}
\label{Lemma modulo m} Let $\mathcal{E}$ and $\mathcal{F}$ be two spectra with
finitely generated homotopy groups in each degree. Let%
\[
f:\mathcal{E\longrightarrow F}%
\]
be a morphism that induces an isomorphism between homotopy groups $\pi
_{i}(-;\,\mathbb{Z}/2^{r})$ for $i\geq n$ and sufficiently large $r$.

\emph{(a) }Then $f$ induces an isomorphism between the $2$-primary torsion of
$\mathcal{E}$ and $\mathcal{F}$ in all degrees $\geq n$.

\emph{(b) }Moreover, if also $f$ induces an isomorphism between rational
homotopy groups for degrees $\geq n$, then it induces an isomorphism between
$\mathcal{E}$ and $\mathcal{F}$ in these degrees after tensoring with
$2$-local integers $\mathbb{Z}_{(2)}$, and hence after tensoring with $2$-adic
integers $\mathbb{Z}_{2}$.
\end{lemma}

\noindent\textbf{Proof.} \textbf{(a) }We consider the map of Bockstein short
exact sequences%
\[%
\begin{array}
[c]{ccccccc}%
0\rightarrow & \pi_{i}(\mathcal{E})/2^{r} & \longrightarrow & \pi
_{i}(\mathcal{E};\,\mathbb{Z}/2^{r}) & \longrightarrow & {}_{2^{r}}\pi
_{i-1}(\mathcal{E}) & \rightarrow0\\
& \downarrow &  & \downarrow^{\cong} &  & \downarrow & \\
0\rightarrow & \pi_{i}(\mathcal{F})/2^{r} & \longrightarrow & \pi
_{i}(\mathcal{F};\,\mathbb{Z}/2^{r}) & \longrightarrow & {}_{2^{r}}\pi
_{i-1}(\mathcal{F}) & \rightarrow0
\end{array}
\]
Evidently, this implies that for all $r$ the map ${}_{2^{r}}\pi_{i-1}%
(\mathcal{E})\rightarrow{}_{2^{r}}\pi_{i-1}(\mathcal{F})$ between $2^{r}%
$-torsion subgroups is surjective whenever $i\geq n$. If now $r$ is
sufficiently large, then injectivity of the map $\pi_{i}(\mathcal{E}%
)/2^{r}\rightarrow\pi_{i}(\mathcal{F})/2^{r}$ yields an isomorphism of the
$2$-primary torsion subgroups, again whenever $i\geq n$.

\noindent\textbf{(b) }Here, following \cite[pg.32]{Sullivan}, we first note
that $\pi_{i}$ commutes with direct limits of coefficient groups, giving an
isomorphism on $\pi_{i}(-;\,\mathbb{Z}_{2^{\infty}})$. Then, from the exact
sequence of homotopy groups associated to the short exact sequence of
coefficients
\[
0\rightarrow\mathbb{Z}_{(2)}\longrightarrow\mathbb{Q}\longrightarrow
\mathbb{Z}_{2^{\infty}}\rightarrow0\text{,}%
\]
we obtain an isomorphism on $\pi_{i}(-;\,\mathbb{Z}_{(2)})$. However, since
$\mathbb{Z}_{(2)}$ is a flat $\mathbb{Z}$-module, the $\mathrm{Tor}$ term
vanishes in the universal coefficient sequence%
\[
0\rightarrow\pi_{i}(\mathcal{E})\otimes_{\mathbb{Z}}\mathbb{Z}_{(2)}%
\longrightarrow\pi_{i}(\mathcal{E};\,\mathbb{Z}_{(2)})\longrightarrow
\mathrm{Tor}_{\mathbb{Z}}(\pi_{i-1}(\mathcal{E}),\,\mathbb{Z}_{(2)}%
)\rightarrow0\text{.}%
\]
Finally, of course, we use the fact that $\mathbb{Z}_{(2)}\otimes_{\mathbb{Z}%
}\mathbb{Z}_{2}\cong\mathbb{Z}_{2}$. \hfill$\Box$\medskip

Following the exposition in \cite[\S 5]{BK}, we start out the proof of Theorem
\ref{theorem1} by choosing an embedding of the field of $q$-adic numbers
$\mathbb{Q}_{q}$ into the complex numbers $\mathbb{C}$ such that the induced
composite map
\[
{}_{\varepsilon}\mathcal{KQ}(\mathbb{Z}_{q})^{c}\rightarrow{}_{\varepsilon
}\mathcal{KQ}(\mathbb{Q}_{q})^{c}\rightarrow{}_{\varepsilon}\mathcal{KQ}%
(\mathbb{C})^{c}%
\]
agrees with Friedlander's Brauer lift ${}_{\varepsilon}\mathcal{KQ}%
(\mathbb{F}_{q})^{c}\rightarrow{}_{\varepsilon}\mathcal{KQ}(\mathbb{C})^{c}$
from \cite{Friedlander} under the rigidity equivalence\textsf{ }between
${}_{\varepsilon}\mathcal{KQ}(\mathbb{Z}_{q})^{c}$ and ${}_{\varepsilon
}\mathcal{KQ}(\mathbb{F}_{q})^{c}$, \textsl{cf}.~\cite[Lemma 5.3]{BK}. This
idea is the hermitian analogue of a widely used construction in algebraic
$K$-theory originating in the works of B{\"{o}}kstedt \cite{Bok},
Dwyer-Friedlander \cite{DF}, and Friedlander \cite{Friedlander}.

The ring maps relating $R_{F}$ to $\mathbb{F}_{q}$, $\mathbb{R}$ and
$\mathbb{C}$ induce the commuting B{\"{o}}kstedt square for hermitian
$K$-theory spectra (\ref{hermitian K Bokstedt square}) in the Introduction.

The same ring maps induce the analogous commuting B{\"{o}}kstedt square for
algebraic $K$-theory spectra (\ref{algebraic K Bokstedt square}) in the
Introduction comprising ${}_{\varepsilon}\mathbb{Z}/2$-equivariant maps.

Denote by ${}_{\varepsilon}\overline{\mathcal{KQ}}(R_{F})$ the homotopy
cartesian product of ${}_{\varepsilon}\mathcal{KQ}(\mathbb{F}_{q})^{c}$ and
$\vee^{r}{}_{\varepsilon}\mathcal{KQ}(\mathbb{R})^{c}$ over $\vee^{r}%
{}_{\varepsilon}\mathcal{KQ}(\mathbb{C})^{c}$, afforded by the B{\"{o}}kstedt
square of hermitian $K$-theory spectra. (${}_{\varepsilon}\overline
{\mathcal{KQ}}(R_{F})$ is thereby connective because of the epimorphism
${}_{\varepsilon}KQ_{0}(\mathbb{R})\twoheadrightarrow{}_{\varepsilon}%
KQ_{0}(\mathbb{C})$. Moreover, since the spectra ${}_{\varepsilon}%
\mathcal{KQ}(\mathbb{F}_{q})^{c}$ and $\vee^{r}{}_{\varepsilon}\mathcal{KQ}%
(\mathbb{R})^{c}$ have finitely generated homotopy groups, so too does
${}_{\varepsilon}\overline{\mathcal{KQ}}(R_{F})$.)

Thus, Theorem \ref{theorem1} becomes the assertion that there is a homotopy
equivalence between the $2$-completed connective spectra associated to the
map
\begin{equation}
{}_{\varepsilon}\mathcal{KQ}(R_{F})\rightarrow{}_{\varepsilon}\overline
{\mathcal{KQ}}(R_{F})\text{.} \label{weakequivalence}%
\end{equation}
We write ${}_{\varepsilon}\overline{KQ}_{n}(R_{F})$ for the homotopy groups of
the target spectrum above.

We prove Theorem \ref{theorem1} by means of the following strategy, previously
adopted in \cite[\S 5]{BK} in the case of the rational field. The
low-dimensional computations in Theorem \ref{lowdimensionalcomputations} below
show that (\ref{weakequivalence}) induces an isomorphism modulo odd torsion,
on integral homotopy groups ${}_{\varepsilon}\pi_{n}$ for $n=0,1$. Moreover,
these homotopy groups are trivial in degree $n=-1$. When this is combined with
the fact that the corresponding algebraic $K$-theory square
(\ref{algebraic K Bokstedt square}) is homotopy cartesian (which is shown in
both \cite{HO} and \cite{Mitchell}, see also Appendix
\ref{section:K-theorybackground}), and the induction methods for hermitian
$K$-groups in \cite[\S 3]{BK}, we deduce that the morphism of spectra
(\ref{weakequivalence}) induces an isomorphism of all homotopy groups both
rationally and with finite $2$-group coefficients. Since both spectra have
their integral homotopy groups finitely generated in all dimensions, we may
apply Lemma \ref{Lemma modulo m} to obtain an isomorphism of $2$-completions
of the integral homotopy groups of these spectra. Then Lemma
\ref{Connective.Spectra}\thinspace(v) finishes the proof.

We must therefore establish the following key computational result which
extends the low-dimensional calculations for the rational integers in
\cite[\S 4]{BK} to every totally real $2$-regular number field.

\begin{theorem}
\label{lowdimensionalcomputations} Let $n=-1,0$ or $1$ and $\varepsilon=\pm1$.
Then the map of homotopy groups
\[
{}_{\varepsilon}\pi_{n}\colon{}_{\varepsilon}KQ_{n}(R_{F})\rightarrow
{}_{\varepsilon}\overline{KQ}_{n}(R_{F})
\]
is an isomorphism, except for $\varepsilon=1$ and $n=0$ where it is an
isomorphism modulo odd torsion.
\end{theorem}

\noindent\textbf{Proof.} The theorem will be proved after many preliminary
lemmas listed below. According to the definition of the spectrum
${}_{\varepsilon}\overline{\mathcal{KQ}}(R_{F})$ as the homotopy cartesian
product of ${}_{\varepsilon}\mathcal{KQ}(\mathbb{F}_{q})^{c}$ and $\vee^{r}%
{}_{\varepsilon}\mathcal{KQ}(\mathbb{R})^{c}$ over $\vee^{r}{}_{\varepsilon
}\mathcal{KQ}(\mathbb{C})^{c}$, there is a naturally induced long exact
sequence
\begin{equation}%
\begin{array}
[c]{ccc}%
\cdots\rightarrow\overset{r}{\bigoplus}{}_{\varepsilon}KQ_{n+1}(\mathbb{C}%
)\rightarrow{}_{\varepsilon}\overline{KQ}_{n}(R_{F})\rightarrow{}%
_{\varepsilon}KQ_{n}(\mathbb{F}_{q})\oplus\overset{r}{\bigoplus}%
{}_{\varepsilon}KQ_{n}(\mathbb{R})\rightarrow\cdots. &  &
\end{array}
\label{longexactsequence1}%
\end{equation}
We shall deal in turn with each of the six cases $\varepsilon=\pm1$,
$n=-1,0,1$, in the remaining lemmas of this section.

\begin{lemma}
\label{lemma:firstcase} The map ${}_{-1}\pi_{1}\colon{}_{-1}KQ_{1}%
(R_{F})\rightarrow{}_{-1}\overline{KQ}_{1}(R_{F})$ is an isomorphism between
trivial groups.
\end{lemma}

\noindent\textbf{Proof.} \label{epsilon=-1, n=1} Recall that ${}_{-1}%
KQ_{1}(A)=0$ when $A=R_{F},\mathbb{F}_{q},\mathbb{R}$ by for example
\cite[Th\'eor\`eme 2.13]{K:AnnScENS} (using \cite[Corollary 12.5]%
{BassMilnorSerre} ). Thus (\ref{longexactsequence1}) shows that it suffices to
note that ${}_{-1}KQ_{2}(\mathbb{C})=\pi_{1}(\mathrm{Sp})=0$ for the infinite
symplectic group $\mathrm{Sp}$.\hfill$\Box$

\begin{lemma}
\label{lemma:secondcase} The map ${}_{-1}\pi_{0}\colon{}_{-1}KQ_{0}%
(R_{F})\rightarrow{}_{-1}\overline{KQ}_{0}(R_{F})\cong\mathbb{Z}$ is an isomorphism.
\end{lemma}

\noindent\textbf{Proof.} \label{epsilon=-1, n=0} By \cite[I (3.5)]{MH} again,
there are isomorphisms ${}_{-1}KQ_{0}(A)\cong\mathbb{Z}$ when $A=R_{F}%
,\mathbb{F}_{q},\mathbb{R}$, and $\mathbb{C}$ detected by the (even) rank of
the free symplectic $A$-inner product space. Now since every ring map
preserves the rank, we obtain a cartesian square:
\begin{equation}%
\begin{array}
[c]{ccc}%
{}_{-1}KQ_{0}(R_{F}) & \rightarrow & \overset{r}{\bigoplus}{}_{-1}%
KQ_{0}(\mathbb{R})\\
\downarrow &  & \downarrow\\
{}_{-1}KQ_{0}(\mathbb{F}_{q}) & \rightarrow & \overset{r}{\bigoplus}{}%
_{-1}KQ_{0}(\mathbb{C})
\end{array}
\label{cartesiandiagram}%
\end{equation}
Combining (\ref{longexactsequence1}), (\ref{cartesiandiagram}), and the fact
that ${}_{-1}KQ_{1}(\mathbb{C})$ is the trivial group, it follows that
${}_{-1}\pi_{0}$ is an isomorphism.\hfill$\Box$

\begin{lemma}
\label{lemma:thirdcase} The map ${}_{1}\pi_{1}\colon{}_{1}KQ_{1}%
(R_{F})\rightarrow{}_{1}\overline{KQ}_{1}(R_{F})$ is an isomorphism.
\end{lemma}

\noindent\textbf{Proof.} \label{epsilon=1, n=1} We first show that the
determinant and the spinor norm of $R_{F}$ induce an isomorphism
\begin{equation}
{}_{1}KQ_{1}(R_{F})\overset{\cong}{\longrightarrow}R_{F}^{\times}%
/(R_{F}^{\times})^{2}\oplus\mathbb{Z}/2. \label{Bassisomorphism}%
\end{equation}
(Moreover, from (\ref{examples for A/A^{2}}) above we note that the right-hand
side reduces to $(\mathbb{Z}/2)^{r+2}$.) To see this, consider the exact
sequence of units, discriminant modules, and Picard groups in \cite[(2.1)]%
{Bass}
\begin{equation}
0\rightarrow\mu_{2}(F)\rightarrow R_{F}^{\times}\overset{(\ )^{2}}%
{\rightarrow}R_{F}^{\times}\rightarrow\mathrm{Discr}(R_{F})\rightarrow
\mathrm{Pic}(R_{F})\overset{2\cdot}{\rightarrow}\mathrm{Pic}(R_{F})\text{.}
\label{exactsequencebass}%
\end{equation}
Combining (\ref{exactsequencebass}) with the $2$-regular assumption on $F$ and
the vanishing of $SK_{1}(R_{F})$ \cite[Corollary 4.3]{BassMilnorSerre}, the
isomorphism in (\ref{Bassisomorphism}) follows from \cite[(4.7.6)]{Bass}. As a
corollary, we deduce that the hyperbolic map $R_{F}^{\times}\cong
K_{1}(F)\rightarrow{}_{1}KQ_{1}(F)\cong R_{F}^{\times}/(R_{F}^{\times}%
)^{2}\oplus\mathbb{Z}/2$ may be identified with the composition of the
quotient map $R_{F}^{\times}\twoheadrightarrow R_{F}^{\times}/(R_{F}^{\times
})^{2}$ by the inclusion of this last group in $R_{F}^{\times}/(R_{F}^{\times
})^{2}\oplus\mathbb{Z}/2$, and so that ${}_{1}W_{1}(R_{F})\cong\mathbb{Z}/2$.
An analogous result holds if we replace $R_{F}$ by $\mathbb{F}_{q},\mathbb{R}$
or $\mathbb{C}$.

The group ${}_{1}KQ_{2}(\mathbb{R})\cong\mathbb{Z}/2\oplus\mathbb{Z}/2$ maps
split surjectively onto ${}_{1}KQ_{2}(\mathbb{C})\cong\mathbb{Z}/2$, as shown
in the first two lemmas of Appendix B. Hence by (\ref{longexactsequence1})
there is a short split exact sequence
\begin{equation}
0\rightarrow{}_{1}\overline{KQ}_{1}(R_{F})\overset{\theta}{\rightarrow}{}%
_{1}KQ_{1}(\mathbb{F}_{q})\oplus\overset{r}{\bigoplus}{}_{1}KQ_{1}%
(\mathbb{R})\overset{\varphi}{\rightarrow}\overset{r}{\bigoplus}{}_{1}%
KQ_{1}(\mathbb{C})\rightarrow0 \label{exactsequence2}%
\end{equation}
where we may choose a splitting of $\varphi$ whose image lies in a direct
summand of $\overset{r}{\bigoplus}{}_{1}KQ_{1}(\mathbb{R})\cong(\mathbb{Z}%
/2)^{2r}$, while $\overset{r}{\bigoplus}{}_{1}KQ_{1}(\mathbb{C})\cong%
(\mathbb{Z}/2)^{r}$. Therefore, we have an induced isomorphism%
\[
{}_{1}\overline{KQ}_{1}(R_{F})\cong{}_{1}KQ_{1}(\mathbb{F}_{q})\oplus
(\mathbb{Z}/2)^{r}\text{.}%
\]
Thus, since ${}_{1}KQ_{1}(\mathbb{F}_{q})\cong\mathbb{Z}/2\oplus\mathbb{Z}/2$,
using (\ref{Bassisomorphism}) and (\ref{exactsequence2}), we deduce that
${}_{1}\overline{KQ}_{1}(R_{F})$ and ${}_{1}KQ_{1}(R_{F})$ are both abstractly
isomorphic to direct sums of $r+2$ copies of $\mathbb{Z}/2$. Therefore, to
finish the proof of the lemma, it suffices to show that ${}_{1}\pi_{1}:{}%
_{1}KQ_{1}(R_{F})\rightarrow{}_{1}\overline{KQ}_{1}(R_{F})$ is surjective. The
argument is now broken into three steps.

\noindent\textbf{Step 1. }For $A=R_{F},\mathbb{F}_{q},\mathbb{R}$ and
$\mathbb{C}$, by Lemma \ref{SK_1(A) = 0}, we have surjections%
\[
{}_{1}KQ_{1}(A)\twoheadrightarrow{}_{1}W_{1}(A)\overset{\cong}%
{\twoheadrightarrow}k_{1}(A)=\{\pm1\}\text{.}%
\]
Let $SKQ_{1}(A)$ denote the kernel of this determinant map ${}_{1}%
KQ_{1}(A)\rightarrow\{\pm1\}$ which is obviously split surjective, so that
${}_{1}KQ_{1}(A)\cong S{}KQ_{1}(A)\oplus\mathbb{Z}/2$. Naturality of the
determinant map implies that the maps between the various ${}_{1}KQ_{1}(A)$
restrict to maps between the corresponding subgroups $SKQ_{1}(A)$. Moreover,
because ${}_{1}W_{1}(A)$ is the cokernel of the hyperbolic map, the map
$K_{1}(A)\rightarrow{}_{1}KQ_{1}(A)$ factors through $SKQ_{1}(A)$.

The computations above now imply that in the commutative diagram
\[%
\begin{array}
[c]{ccc}%
R_{F}^{\times}\cong K_{1}(R_{F}) & \longrightarrow & K_{1}(\mathbb{F}%
_{q})\oplus\overset{r}{\bigoplus}K_{1}(\mathbb{R})\\
\downarrow &  & \downarrow\\
R_{F}^{\times}/(R_{F}^{\times})^{2}\cong SKQ_{1}(R_{F}) & \longrightarrow &
SKQ_{1}(\mathbb{F}_{q})\oplus\overset{r}{\bigoplus}SKQ_{1}(\mathbb{R})
\end{array}
\]
both the vertical and top horizontal arrows are surjective. Therefore, the
lower horizontal map is also surjective.

\noindent\textbf{Step 2. }In a symmetric way, we introduce a subgroup
$S\overline{KQ}_{1}(R_{F})$ of ${}_{1}\overline{KQ}_{1}(R_{F})$ as the kernel
of the composition ${}_{1}\overline{KQ}_{1}(R_{F})\cong{}_{1}KQ_{1}%
(\mathbb{F}_{q})\oplus(\mathbb{Z}/2)^{r}\longrightarrow{}_{1}KQ_{1}%
(\mathbb{F}_{q})\overset{\mathrm{det}}{\longrightarrow}\{\pm1\}$ which is a
split surjection. Therefore, we have a decomposition of both ${}_{1}%
KQ_{1}(R_{F})$ and ${}_{1}\overline{KQ}_{1}(R_{F})$ as compatible direct sums:%
\[
{}_{1}KQ_{1}(R_{F})\cong SKQ_{1}(R_{F})\oplus\mathbb{Z}/2\text{\quad and\quad
}{}_{1}\overline{KQ}_{1}(R_{F})\cong S{}\overline{KQ}_{1}(R_{F})\oplus
\mathbb{Z}/2\text{,}%
\]
and it is enough to prove the surjectivity of the induced map $\sigma:{}%
_{1}SKQ_{1}(R_{F})\rightarrow{}_{1}\overline{SKQ}_{1}(R_{F})$.

\noindent\textbf{Step 3. }Finally, consider the commuting diagram:
\[%
\begin{array}
[c]{ccccc}%
SKQ_{1}(R_{F}) & \overset{\sigma}{\longrightarrow} & S\overline{KQ}_{1}%
(R_{F}) & \overset{\gamma}{\longrightarrow} & SKQ_{1}(\mathbb{F}_{q}%
)\oplus\overset{r}{\bigoplus}SKQ_{1}(\mathbb{R})\\
\downarrow &  & \downarrow &  & \downarrow\\
{}_{1}KQ_{1}(R_{F}) & \overset{{}_{1}\pi_{1}}{\longrightarrow} & {}%
_{1}\overline{KQ}_{1}(R_{F}) & \overset{\theta}{\longrightarrow} & {}%
_{1}KQ_{1}(\mathbb{F}_{q})\oplus\overset{r}{\bigoplus}{}_{1}KQ_{1}(\mathbb{R})
\end{array}
\]
From (\ref{exactsequence2}), $\theta$ is a monomorphism; so, its restriction
$\gamma$ is also injective. Meanwhile, from Step 1, the composite $\gamma
\circ\sigma$ is surjective. It follows that $\sigma$ is surjective too, as we
sought.\hfill$\Box$

\begin{lemma}
\label{lemma:fourthcase} The map ${}_{1}\pi_{0}\colon{}_{1}KQ_{0}%
(R_{F})\rightarrow{}_{1}\overline{KQ}_{0}(R_{F})$ is an isomorphism modulo odd torsion.
\end{lemma}

\noindent\textbf{Proof.} Since $\mathrm{Pic}(A)$ is an odd torsion group we
obtain for $A=R_{F},\mathbb{F}_{q},\mathbb{R}$, and $\mathbb{C}$ the short
exact sequence (modulo odd torsion)
\begin{equation}
0\rightarrow W^{\prime}(A)\rightarrow{}_{1}KQ_{0}(A)\rightarrow K_{0}%
(A)\rightarrow0.\label{ses W' -> KQ -> K}%
\end{equation}
From the last three of these four cases we obtain the vertical map of short
exact sequences:
\[%
\begin{array}
[c]{ccccccc}%
0\rightarrow & W^{\prime}(\mathbb{F}_{q})\oplus\bigoplus\limits^{r}W^{\prime
}(\mathbb{R}) & \rightarrow & {}_{1}KQ_{0}(\mathbb{F}_{q})\oplus
\bigoplus\limits^{r}{}_{1}KQ_{0}(\mathbb{R}) & \rightarrow & K_{0}%
(\mathbb{F}_{q})\oplus\bigoplus\limits^{r}K_{0}(\mathbb{R}) & \rightarrow0\\
& \downarrow &  & \downarrow &  & \downarrow & \\
0\rightarrow & \bigoplus\limits^{r}W^{\prime}(\mathbb{C}) & \rightarrow &
\bigoplus\limits^{r}{}_{1}KQ_{0}(\mathbb{C}) & \rightarrow & \bigoplus
\limits^{r}K_{0}(\mathbb{C}) & \rightarrow0
\end{array}
\]
By \cite{HO} and \cite{Mitchell}, $K_{0}(R_{F})$ is the kernel of the
rightmost vertical map. Further, the bottom right horizontal epimorphism maps
between two copies of $\mathbb{Z}^{r}$, making the coWitt group $W^{\prime
}(\mathbb{C})$ trivial. As already noted, ${}_{1}KQ_{1}(\mathbb{R})$ maps onto
${}_{1}KQ_{1}(\mathbb{C})$, making ${}_{1}\overline{KQ}_{0}(R_{F})$ the kernel
of the center vertical map. Therefore, the desired isomorphism between ${}%
_{1}KQ_{0}(R_{F})$ and ${}_{1}\overline{KQ}_{0}(R_{F})$ follows by combining
the short exact sequence (\ref{ses W' -> KQ -> K}) for $A=R_{F}$ and
Proposition \ref{proposition:keyinput}.$\hfill\Box\smallskip$

Before dealing with the remaining cases $n=-1$ with $\varepsilon=\pm1$, we
need a proposition which is interesting by itself.

\begin{proposition}
\label{Dedekind forgetful}Let $A$ be a Dedekind ring and $\varepsilon=\pm1$
and let us consider the forgetful map
\[
{}_{\varepsilon}KQ_{0}(A)\overset{\varphi}{\longrightarrow}K_{0}%
(A)=\mathbb{Z}\oplus\mathrm{Pic}(A)\text{.}%
\]
Its image on the first summand is $\mathbb{Z}$ if $\varepsilon=1$ and
$2\mathbb{Z}$ if $\varepsilon=-1$. Its image on the second summand is a group
of order at most $2$. Therefore, if $A$ is a $2$-regular ring of integers, the
image of $\varphi$ lies in the first summand.
\end{proposition}

\noindent\textbf{Proof.} The projection to the $\mathbb{Z}$ summand is given
by the rank. On the other hand, if $E$ is equipped with a symmetric or
antisymmetric form, it is isomorphic to its dual. The duality on
$\mathrm{Pic}(A)$ is given by $I\longmapsto I^{\ast}=I^{-1}$. This implies
that the image of $E$ in $\mathrm{Pic}(A)$ is at most of order $2$. In
particular, when $\mathrm{Pic}(A)$ is odd torsion, the image of $\varphi$ is
included in the $\mathbb{Z}$ summand. \hfill$\Box$

For the cases $n=-1$ with $\varepsilon=\pm1$, we are going to show that
${}_{\varepsilon}KQ_{-1}(R_{F})$ is trivial, as is ${}_{\varepsilon}%
\overline{KQ}_{-1}(R_{F})$ because the spectrum ${}_{\varepsilon}%
\overline{\mathcal{KQ}}(R_{F})$ is connective. As in \cite{BK}, we exploit the
two exact sequences
\[
K_{1}(R_{F})\longrightarrow{}_{-\varepsilon}KQ_{1}(R_{F})\longrightarrow
{}_{-\varepsilon}U_{0}(R_{F})\longrightarrow K_{0}(R_{F})\longrightarrow
{}_{-\varepsilon}KQ_{0}(R_{F}),
\]
and (via the fundamental theorem proved in \cite{K:AnnM112fun})%
\[
{}_{\varepsilon}KQ_{0}(R_{F})\overset{\varphi}{\longrightarrow}K_{0}%
(R_{F})\longrightarrow{}_{-\varepsilon}U_{0}(R_{F})\longrightarrow
{}_{\varepsilon}KQ_{-1}(R_{F})\longrightarrow K_{-1}(R_{F})=0\text{.}%
\]

\begin{lemma}
\label{lemma:epsilon=1, n=-1} The map ${}_{1}\pi_{-1}\colon{}_{1}KQ_{-1}%
(R_{F})\rightarrow{}_{1}\overline{KQ}_{-1}(R_{F})$ is an isomorphism between
trivial groups.
\end{lemma}

\noindent\textbf{Proof.} For $\varepsilon=1$, ${}_{-\varepsilon}KQ_{1}%
(R_{F})=0$ by a well-known theorem of Bass, Milnor and Serre \cite[Corollary
12.5]{BassMilnorSerre}; and the map $K_{0}(R_{F})\rightarrow{}_{-\varepsilon
}KQ_{0}(R_{F})$ has its kernel identified with $\mathrm{Pic}(R_{F})$ since
${}_{-1}KQ_{0}(R_{F})=\mathbb{Z}$. Hence ${}_{-\varepsilon}U_{0}(R_{F}%
)\cong\mathrm{Pic}(R_{F})$. On the other hand, according to the previous
proposition, the cokernel of $\varphi$ is also identified with $\mathrm{Pic}%
(R_{F})$. The second exact sequence therefore provides an injective map from
$\mathrm{Pic}(R_{F})$ to itself. Since $\mathrm{Pic}(R_{F})$ is finite, this
map is bijective, and it follows that ${}_{1}KQ_{-1}(R_{F})$ is the trivial
group.\hfill$\Box$

\begin{lemma}
\label{lemma:epsilon=-1, n=-1} The map ${}_{-1}\pi_{-1}\colon{}_{-1}%
KQ_{-1}(R_{F})\rightarrow{}_{-1}\overline{KQ}_{-1}(R_{F})$ is an isomorphism
between trivial groups\textsf{.}
\end{lemma}

\noindent\textbf{Proof.} In order to compute ${}_{-1}KQ_{-1}(R_{F}),$ let us
first work modulo odd torsion, which makes the map $K_{0}(R_{F}%
)\longrightarrow{}_{-\varepsilon}KQ_{0}(R_{F})$ injective. The first exact
sequence above written for $\varepsilon=-1$ shows that ${}_{1}U_{0}%
(R_{F})=\mathbb{Z}/2$, because ${}_{1}W_{1}(R_{F})=\mathbb{Z}/2$ by
(\ref{W1}). On the other hand, the cokernel of the map $\varphi$ in the second
exact sequence is also $\mathbb{Z}/2$ by Proposition \ref{Dedekind forgetful}.
It follows that $_{-1}KQ_{-1}(R_{F})=0$, and so the map ${}_{-1}KQ_{-1}%
(R_{F})\rightarrow{}_{-1}\overline{KQ}_{-1}(R_{F})$ is a homomorphism between
odd order finite groups.

Finally, according to Proposition \ref{Wodd}, the odd torsion of the Witt
groups in degree $-1$ is trivial. Since $K_{-1}(R_{F})=0$, it follows that the
odd torsion of ${}_{-1}KQ_{-1}(R_{F})$ is also trivial.\hfill$\Box\smallskip$

These various lemmas conclude the proof of Theorem
\ref{lowdimensionalcomputations}, and hence of Theorem \ref{theorem1} with the
exception of the following proposition which was announced at the beginning:

\begin{proposition}
\label{Finiteness}Let $F$ be a number field and let $R$ be a ring of
$S$-integers in $F$, where $S$ is a finite set of primes containing the dyadic
ones. Then the groups ${}_{\varepsilon}KQ_{n}(R)$ are finitely generated for
$\varepsilon=\pm1$ and $n\in\mathbb{Z}$.
\end{proposition}

\noindent\textbf{Proof.} The same arguments used in the proofs of Lemmas
\ref{lemma:firstcase} and \ref{lemma:secondcase} show that ${}_{-1}%
KQ_{0}(R)\equiv\mathbb{Z}$ and $_{-1}KQ_{1}(R)=0$. According to \cite[(2.1)]%
{Bass}, $_{1}KQ_{1}(R)$ ($=KO_{1}(R)$ in Bass's notation) is inserted in an
exact sequence between two finitely generated groups $K\mathrm{SL}_{1}(R)$ and
$\mathrm{Disc}(R)\oplus\mathbb{Z}/2$, where $\mathrm{Disc}(R)$ is described in
\cite[pg. 156.]{Bass}. Therefore, $_{1}KQ_{1}(R)$ and {}$_{1}W_{1}(R)$ are
finitely generated.

On the other hand, it is well known (see \textsl{e.g.} \cite[ Corollary 3.3,
pg. 93]{MH}), that the canonical map between classical Witt groups
$W(R)\longrightarrow W(F)$ is injective. Moreover, any element of $W(F)$ is
determined by the classical invariants which are the rank, signature,
discriminant and Hasse-Witt invariant. Since $S$ is finite, the discriminant
and the Hasse-Witt invariant computed on the elements of $W(R)$, considered as
a subgroup of $W(F)$, can take only a finite set of values. Moreover, the
signature takes integral values defined by the various real embeddings$.$
Therefore, the group $W(R)$ is finitely generated.

In order to deal with the higher Witt groups ${}_{\varepsilon}W_{n}(R)$ and
coWitt groups ${}_{\varepsilon}W_{n}^{^{\prime}}(R)$ for $n\in\mathbb{Z}$, we
use the following two basic tools:

1) According to Quillen \cite{Quillenfinite} the groups $K_{n}(R)$ are
finitely generated. Therefore, the Tate cohomology groups $k_{n}(R)$ and
$k_{n}^{^{\prime}}(R)$ are also finitely generated.

2) There is a $12$-term exact sequence detailed in \cite[pg. 278.]%
{K:AnnM112fun} which shows by double induction from the cases $\varepsilon
=\pm1$ and $n=0,1$, that the groups ${}_{\varepsilon}W_{n}(R)$ and
${}_{\varepsilon}W_{n}^{^{\prime}}(R)$ are finitely generated for all values
of $n\in\mathbb{Z}$ and $\varepsilon=\pm1$.

Finally, from 1) and 2), we deduce that ${}_{\varepsilon}KQ_{n}(R)$, which
lies in an exact sequence between $K_{n}(R)$ and ${}_{\varepsilon}W_{n}(R)$,
is finitely generated.\textbf{\hfill}$\Box$\medskip

\begin{remark}
\label{localization}1. By means of Lemma \ref{Lemma modulo m}\thinspace(b),
the proof of Theorem \ref{theorem1} also works when dealing with
$2$-localizations instead of $2$-completions. Likewise, Theorems
\ref{theorem2} and \ref{theorem4} are also true in the framework of
$2$-localizations. However, Theorem \ref{theorem3} is only true for completions.

2. In this paper, we work with the spectra ${}\mathcal{K}(\mathbb{R})_{\#}%
^{c}$ and ${}\mathcal{K}(\mathbb{C})_{\#}^{c}$ and their hermitian analogs.
The Suslin equivalence \cite{Suslin :local fields} with ${}\mathcal{K}%
(\mathbb{R}^{\delta})_{\#}$ and ${}\mathcal{K}(\mathbb{C}^{\delta})_{\#}$
respectively, where $\delta$ means the discrete topology, enables these
spectra to be replaced by those of the discrete rings when considering
completions. However, this method fails if we deal with localizations instead
of completions.
\end{remark}

\section{Splitting results}

\label{section:splittingresults} The purpose of this section is to prove some
generic splitting results employed in the proofs of Theorems \ref{theorem2}
and \ref{theorem4}. \vspace{0.1in}

To start with, fix some residue field $\mathbb{F}_{q}$ of $R_{F}$ as in
paragraph prior to the statement of Theorem \ref{theorem1} in the
Introduction, and define the spectrum ${}_{\varepsilon}\mathcal{KQ}%
(\overline{R}_{F})$ (this is just a convenient notation, which we shall use
below for $K$-theory and $V$-theory too) by the homotopy cartesian square:
\[
\xymatrix{ {}_\varepsilon\mathcal{KQ}(\overline{R}_{F})\ar@{}[d]|-{\hbox
{\large{$\downarrow$}}} \ar@{}[r]|-{\hbox{\large{$\rightarrow$}}}&
{}_\varepsilon\mathcal{KQ}(\mathbb{R})^c\ar@{}[d]|-{\hbox{\large{$\downarrow
$}}}\\ {}_\varepsilon\mathcal{KQ}(\mathbb{F}_q)^c\ar@{}[r]|-{\hbox
{\large{$\rightarrow$}}} & {}_\varepsilon\mathcal{KQ}(\mathbb{C})^c }\label{RFoverbar}%
\]
It is connective because of the epimorphism ${}_{\varepsilon}KQ_{0}(%
\mathbb{R}%
)\rightarrow{}_{\varepsilon}KQ_{0}(%
\mathbb{C}%
)$. Similarly, define $\mathcal{K}(\overline{R}_{F})$ by the same type of
homotopy cartesian square by replacing ${}_{\varepsilon}\mathcal{KQ}$ by
$\mathcal{K}.$ It is also connective because of the epimorphism $K_{0}(%
\mathbb{R}%
){}\rightarrow K_{0}(%
\mathbb{C}%
)$. Finally, we define ${}_{\varepsilon}\mathcal{V}(\overline{R}_{F})$ by the
similar homotopy cartesian square replacing ${}_{\varepsilon}\mathcal{KQ}$ by
${}_{\varepsilon}\mathcal{V}$. It too is connective, because of the
epimorphism ${}_{\varepsilon}V_{0}(%
\mathbb{R}%
)\rightarrow{}_{\varepsilon}V_{0}(%
\mathbb{C}%
)$ proved in Lemmas \ref{_1VR}, \ref{_1VC} and \ref{_-1VRC} of Appendix B.

The first homotopy cartesian square can be recast as a homotopy fiber
sequence${}$
\begin{equation}
_{\varepsilon}\mathcal{KQ}(\overline{R}_{F})\longrightarrow{}_{\varepsilon
}\mathcal{KQ}(%
\mathbb{R}
)^{c}{}\vee{}_{\varepsilon}\mathcal{KQ}(\mathbb{F}_{q})^{c}\longrightarrow
{}_{\varepsilon}\mathcal{KQ}(%
\mathbb{C}
)^{c}\text{,} \label{_epsilonKQ(R^-_F) fibration}%
\end{equation}
and similarly for the others.

In the next step we form the naturally induced diagram with horizontal
homotopy fiber sequences:
\begin{equation}%
\begin{array}
[c]{ccccc}%
_{\varepsilon}\mathcal{KQ}(\overline{R}_{F}) & \longrightarrow &
{}_{\varepsilon}\mathcal{KQ}(%
\mathbb{R}
)^{c}{}\vee{}_{\varepsilon}\mathcal{KQ}(\mathbb{F}_{q})^{c} & \longrightarrow
& _{\varepsilon}\mathcal{KQ}(%
\mathbb{C}
)^{c}\\
\downarrow &  & \downarrow^{\nabla\vee\mathrm{id}} &  & \downarrow^{\nabla}\\
_{\varepsilon}\overline{\mathcal{KQ}}(R_{F}) & \longrightarrow &
{\displaystyle\bigvee\nolimits^{r}}
{}_{\varepsilon}\mathcal{KQ}(%
\mathbb{R}
)^{c}{}\vee{}_{\varepsilon}\mathcal{KQ}(\mathbb{F}_{q})^{c} & \longrightarrow
&
{\displaystyle\bigvee\nolimits_{\varepsilon}^{r}}
\mathcal{KQ}(%
\mathbb{C}
)^{c}\\
\downarrow &  & \downarrow &  & \downarrow\\%
{\displaystyle\bigvee\nolimits^{r-1}}
{}_{\varepsilon}\mathcal{F} & \longrightarrow &
{\displaystyle\bigvee\nolimits^{r-1}}
{}_{\varepsilon}\mathcal{KQ}(%
\mathbb{R}
)^{c}{} & \longrightarrow &
{\displaystyle\bigvee\nolimits_{\varepsilon}^{r-1}}
\mathcal{KQ}(%
\mathbb{C}
)^{c}%
\end{array}
\label{3 x 3}%
\end{equation}
Note that ${}_{\varepsilon}\mathcal{F}$ is connective, again from the
epimorphism ${}_{\varepsilon}KQ_{0}(%
\mathbb{R}
)\rightarrow{}_{\varepsilon}KQ_{0}(%
\mathbb{C}
)$. Since the two right-hand columns are also homotopy fiber sequences, the
same holds for the left-hand column, \textsl{cf.}~\cite[Lemma 2.1]{CMN}. Using
the compatibility of the two evident splittings of the right-hand columns, we
obtain a splitting of ${}_{\varepsilon}\mathcal{KQ}(\overline{R}%
_{F})\rightarrow{}_{\varepsilon}\overline{\mathcal{KQ}}(R_{F})$. In other
words,%
\[
{}_{\varepsilon}\overline{\mathcal{KQ}}(R_{F})\simeq{}_{\varepsilon
}\mathcal{KQ}(\overline{R}_{F})\vee\bigvee^{r-1}{}_{\varepsilon}%
\mathcal{F}\text{.}%
\]

Let us define the connective spectrum $\overline{\mathcal{K}}(R_{F})$ as the
homotopy pull-back of the diagram%
\[%
\begin{array}
[c]{ccc}%
\overline{\mathcal{K}}(R_{F}) & \rightarrow & {\displaystyle\bigvee
\limits^{r}}{}\mathcal{K}(\mathbb{R})^{c}\\
\downarrow &  & \downarrow\\
\mathcal{K}(\mathbb{F}_{q}) & \rightarrow & {\displaystyle\bigvee\limits^{r}%
}\mathcal{K}(\mathbb{C})^{c}%
\end{array}
\]
It is connective because the map $K_{0}(%
\mathbb{R}
)\longrightarrow K_{0}(%
\mathbb{C}
)$ is surjective. By an argument similar to before, we can split
$\overline{\mathcal{K}}(R_{F})$ as $\overline{\mathcal{K}}(R_{F}%
)\simeq\mathcal{K}(\overline{R}_{F})\vee\bigvee^{r-1}\Omega^{-1}%
\mathcal{K}(\mathbb{R})^{c}$, where $\Omega^{-1}\mathcal{K}(\mathbb{R}%
)^{c}=(\Omega^{-1}\mathcal{K}(\mathbb{R}))^{c}$ (since $K_{-1}(%
\mathbb{R}
)=0$) is the homotopy fiber of the map $\mathcal{K}(\mathbb{R})^{c}%
\rightarrow\mathcal{K}(\mathbb{C})^{c}$ by a direct consequence of a
well-known result due to Bott: see for instance \cite[Section III.5]%
{K:ktheorybook}. This splitting will be used incidentally in the computation
of one $V$-group later on. It also implies the explicit computation of the
groups $K_{n}(R_{F})$ listed in Theorem $1.3.$

In ${}_{1}\mathcal{KQ}$-theory, we make use of the following lemma, which is a
consequence of (\ref{3 x 3}) and the homotopy equivalences $_{1}\mathcal{KQ}(%
\mathbb{R}
)\simeq\mathcal{K}(%
\mathbb{R}
)\vee\mathcal{K}(%
\mathbb{R}
)$ and $_{1}\mathcal{KQ}(%
\mathbb{C}
)\simeq\mathcal{K}(%
\mathbb{R}
)$ proved in Appendix B which give the homotopy type of $_{1}\mathcal{F}$ as
$\mathcal{K}(\mathbb{R})^{c}$.

\begin{lemma}
\label{_1KQ(R_F) splits}There is a homotopy equivalence
\[
{}_{1}\overline{\mathcal{KQ}}(R_{F})\simeq{}_{1}\mathcal{KQ}(\overline{R}%
_{F})\vee%
{\displaystyle\bigvee\nolimits^{r-1}}
\mathcal{K}(\mathbb{R})^{c}.\vspace{-20pt}%
\]
\hfill$\Box\smallskip$
\end{lemma}

Next, we consider the fiber ${}_{-1}\mathcal{F}$ of
\[
{}_{-1}\mathcal{KQ}(\mathbb{R})^{c}=\mathcal{K}(\mathbb{C})^{c}\longrightarrow
{}_{-1}\mathcal{KQ}(\mathbb{C})^{c}=\mathcal{K}(\mathbb{H})^{c}%
\]
(where $\mathbb{H}$ refers to the quaternions with the usual topology). If we
write the commutative diagram due to Bott%
\[%
\begin{array}
[c]{ccc}%
\mathcal{K}(\mathbb{C})^{c} & \longrightarrow & \mathcal{K}(\mathbb{H})^{c}\\
\downarrow^{\simeq} &  & \downarrow^{\simeq}\\
(\Omega^{4}(\mathcal{K}(\mathbb{C})))^{c} & \longrightarrow & (\Omega
^{4}(\mathcal{K}(\mathbb{R})))^{c}%
\end{array}
\]
we see that the homotopy fiber of the map $_{-1}\mathcal{KQ}(\mathbb{R}%
)^{c}\rightarrow{}_{-1}\mathcal{KQ}(\mathbb{C})^{c}$ may be identified with
the homotopy fiber $(\Omega^{6}\mathcal{K}(\mathbb{R}))^{c}$ of $(\Omega
^{4}(\mathcal{K}(\mathbb{C})))^{c}\rightarrow(\Omega^{4}(\mathcal{K}%
(\mathbb{R})))^{c}$.

\begin{lemma}
\label{_-1KQ(R_F) splits}There is a homotopy equivalence
\[
{}_{-1}\overline{\mathcal{KQ}}(R_{F})\simeq{}_{-1}\mathcal{KQ}(\overline
{R}_{F})\vee\bigvee^{r-1}(\Omega^{6}\mathcal{K}(\mathbb{R}))^{c}%
.\vspace{-20pt}%
\]
\hfill$\Box\smallskip$
\end{lemma}

The next proposition is a consequence of the two previous lemmas and Theorem
\ref{theorem1}.

\begin{proposition}
There are homotopy equivalences of $2$-completed connective spectra%
\[
{}_{1}\mathcal{KQ}(R_{F})_{\#}^{c}\simeq{}_{1}\mathcal{KQ}(\overline{R}%
_{F})_{\#}\vee%
{\displaystyle\bigvee\nolimits^{r-1}}
\mathcal{K}(\mathbb{R})_{\#}^{c}%
\]
and%
\[
{}_{-1}\mathcal{KQ}(R_{F})_{\#}^{c}\simeq{}_{-1}\mathcal{KQ}(\overline{R}%
_{F})_{\#}\vee\bigvee^{r-1}(\Omega^{6}\mathcal{K}(\mathbb{R}))_{\#}%
^{c}.\vspace{-20pt}%
\]
$\hfill\Box\smallskip$
\end{proposition}

In order to obtain the diagram in ${}_{\varepsilon}\mathcal{V}$-theory
corresponding to (\ref{3 x 3}), we start with the diagram
\begin{equation}%
\begin{array}
[c]{ccccc}%
{}_{\varepsilon}\mathcal{V}(R_{F})^{c} & \longrightarrow & \bigvee^{r}%
{}_{\varepsilon}\mathcal{V}(\mathbb{R})^{c}\vee{}_{\varepsilon}\mathcal{V}%
(\mathbb{F}_{q})^{c} & \longrightarrow & \bigvee^{r}{}_{\varepsilon
}\mathcal{V}(\mathbb{C})^{c}\\
\downarrow &  & \downarrow &  & \downarrow\\
{}_{\varepsilon}\mathcal{KQ}(R_{F})^{c} & \longrightarrow & \bigvee^{r}%
{}_{\varepsilon}\mathcal{KQ}(\mathbb{R})^{c}\vee{}_{\varepsilon}%
\mathcal{KQ}(\mathbb{F}_{q})^{c} & \longrightarrow & \bigvee^{r}%
{}_{\varepsilon}\mathcal{KQ}(\mathbb{C})^{c}\\
\downarrow &  & \downarrow &  & \downarrow\\
\mathcal{K}(R_{F}) & \longrightarrow & \bigvee^{r}\mathcal{K}(\mathbb{R}%
)^{c}\vee\mathcal{K}(\mathbb{F}_{q}) & \longrightarrow & \bigvee
^{r}\mathcal{K}(\mathbb{C})^{c}%
\end{array}
\label{3 x 3 V}%
\end{equation}
in which by definition all three columns are fiber sequences; by Theorem
\ref{theorem1} and its counterpart in algebraic $K$-theory
(\ref{algebraic K Bokstedt square}). The lower two rows are fiber sequences,
and hence the top row is as well \cite[Lemma 2.1]{CMN}. A similar argument
reveals that ${}_{\varepsilon}\mathcal{V}(\overline{R}_{F})$ is the homotopy
fiber of the map ${}_{\varepsilon}\mathcal{KQ}(\overline{R}_{F})\rightarrow
\mathcal{K}(\overline{R}_{F})$.

Let us define the connective spectrum ${}_{\varepsilon}\overline{\mathcal{V}%
}(R_{F})$ as the homotopy pull-back of the diagram%
\[%
\begin{array}
[c]{ccc}%
{}_{\varepsilon}\overline{\mathcal{V}}(R_{F}) & \rightarrow &
{\displaystyle \bigvee\limits^{r}}{}{}_{\varepsilon}\mathcal{V}(\mathbb{R}%
)^{c}\\
\downarrow &  & \downarrow\\
{}_{\varepsilon}\mathcal{V}(\mathbb{F}_{q}) & \rightarrow &
{\displaystyle \bigvee\limits^{r}}{}_{\varepsilon}\mathcal{V}(\mathbb{C})^{c}%
\end{array}
\]
This spectrum is connective because the map ${}_{\varepsilon}V_{0}(%
\mathbb{R}%
)\longrightarrow$ ${}_{\varepsilon}V_{0}(%
\mathbb{C}%
)=0$ is surjective according to Lemmas \ref{_1VC} and \ref{_-1VRC} of Appendix B.

By considering the homotopy fibers of the map from Diagram \ref{3 x 3} to its
$\mathcal{K}$-counterpart, we now obtain the commuting diagram
\begin{equation}%
\begin{array}
[c]{ccccc}%
_{\varepsilon}\mathcal{V}(\overline{R}_{F}) & \longrightarrow & {}%
_{\varepsilon}\mathcal{V}(%
\mathbb{R}
)^{c}{}\vee{}_{\varepsilon}\mathcal{V}(\mathbb{F}_{q})^{c} & \longrightarrow &
_{\varepsilon}\mathcal{V}(%
\mathbb{C}
)^{c}\\
\downarrow &  & \downarrow^{\nabla\vee\mathrm{id}} &  & \downarrow^{\nabla}\\
_{\varepsilon}\overline{\mathcal{V}}(R_{F}) & \longrightarrow &
{\displaystyle\bigvee\nolimits^{r}}
{}_{\varepsilon}\mathcal{V}(%
\mathbb{R}
)^{c}{}\vee{}_{\varepsilon}\mathcal{V}(\mathbb{F}_{q})^{c} & \longrightarrow
&
{\displaystyle\bigvee\nolimits_{\varepsilon}^{r}}
\mathcal{V}(%
\mathbb{C}
)^{c}\\
\downarrow &  & \downarrow &  & \downarrow\\%
{\displaystyle\bigvee\nolimits^{r-1}}
{}_{\varepsilon}\mathcal{G} & \longrightarrow &
{\displaystyle\bigvee\nolimits^{r-1}}
{}_{\varepsilon}\mathcal{V}(%
\mathbb{R}
)^{c}{} & \longrightarrow &
{\displaystyle\bigvee\nolimits_{\varepsilon}^{r-1}}
\mathcal{V}(%
\mathbb{C}
)^{c}%
\end{array}
\label{3 x 3 V 2nd}%
\end{equation}
wherein we conclude that the first column is a homotopy fiber sequence, since
all other columns and rows are. As for (\ref{3 x 3}), we deduce that the first
column is split. Thus, a complete determination of ${}_{\varepsilon}%
\overline{\mathcal{V}}(R_{F})$ in terms of ${}_{\varepsilon}\mathcal{V}%
(\overline{R}_{F})$ requires only a computation of the summand {}%
$_{\varepsilon}\mathcal{G}$ which is (by definition) the homotopy fiber of
${}_{\varepsilon}\sigma:{}_{\varepsilon}\mathcal{V}(\mathbb{R})^{c}%
\rightarrow{}_{\varepsilon}\mathcal{V}(\mathbb{C})^{c}$. The following is
included in Appendix B, as Lemmas \ref{_1VR -> C} and \ref{_-1V(R) -> _-1V(C)}.

\begin{lemma}
\label{Fiber V(R) to V(C)}The homotopy fiber {}$_{\varepsilon}\mathcal{G}$ of
${}_{\varepsilon}\sigma:{}_{\varepsilon}\mathcal{V}(\mathbb{R})^{c}%
\rightarrow{}_{\varepsilon}\mathcal{V}(\mathbb{C})^{c}$ is
\[
\left\{
\begin{array}
[c]{lll}%
\mathcal{K}(\mathbb{R})^{c}\vee\mathcal{K}(\mathbb{R})^{c} & \quad &
\varepsilon=1,\\
\mathcal{K}(\mathbb{C})^{c} &  & \varepsilon=-1.
\end{array}
\right.  \vspace{-20pt}%
\]
\hfill$\Box\smallskip$
\end{lemma}

The above computations therefore imply the following two lemmas.

\begin{lemma}
\label{_1V(R_F) splits} There is a homotopy equivalence
\[
{}_{1}\overline{\mathcal{V}}(R_{F})\simeq{}_{1}\mathcal{V}(\overline{R}%
_{F})\vee\bigvee^{r-1}\mathcal{K}(\mathbb{R})^{c}\vee\bigvee^{r-1}%
\mathcal{K}(\mathbb{R})^{c}.\vspace{-20pt}%
\]
\hfill$\Box\smallskip$
\end{lemma}

\begin{lemma}
\label{_-1V(R_F) splits} There is a homotopy equivalence
\[
{}_{-1}\overline{\mathcal{V}}(R_{F})\simeq{}_{-1}\mathcal{V}(\overline{R}%
_{F})\vee\bigvee^{r-1}\mathcal{K}(\mathbb{C})^{c}\text{.}\vspace{-20pt}%
\]
\hfill$\Box\smallskip$
\end{lemma}

Let us now consider the following diagram of fibrations%
\[%
\begin{array}
[c]{ccccc}%
{}_{\varepsilon}\mathcal{V}(R_{F})_{\#}^{c} & \longrightarrow & {}%
_{\varepsilon}\mathcal{KQ}(R_{F})_{\#}^{c} & \longrightarrow & \mathcal{K}%
(R_{F})_{\#}^{c}\\
\downarrow^{\alpha} &  & \downarrow^{\beta} &  & \downarrow^{\gamma}\\
{}_{\varepsilon}\overline{\mathcal{V}}(R_{F})_{\#} & \longrightarrow &
{}_{\varepsilon}\overline{\mathcal{KQ}}(R_{F})_{\#} & \longrightarrow &
\overline{\mathcal{K}}(R_{F})_{\#}%
\end{array}
\]
Since $\beta$ and $\gamma$ are homotopy equivalences, $\alpha$ is a homotopy
equivalence. To finish the computations of the $KQ$ and $V$-groups, it remains
to compute explicitly the groups ${}_{\varepsilon}KQ_{n}(\overline{R}_{F}%
)=\pi_{n}({}_{\varepsilon}\mathcal{KQ}(\overline{R}_{F}))$ and ${}%
_{\varepsilon}V_{n}(\overline{R}_{F})=\pi_{n}({}_{\varepsilon}\mathcal{V}%
(\overline{R}_{F}))$ in the next two sections.

\section{Proof of Theorem \ref{theorem2}}

\label{section:proofoftheorem2} The splitting results Lemma
\ref{_1KQ(R_F) splits} and Lemma \ref{_-1KQ(R_F) splits} in
\S \ref{section:splittingresults} show that in order to prove Theorem
\ref{theorem2}, it suffices to compute the groups ${}_{\varepsilon}%
KQ_{n}(\overline{R}_{F})=\pi_{n}(_{\varepsilon}\mathcal{KQ}(\overline{R}%
_{F}))$, and then sum with $r-1$ copies of the well-known $K$-groups of
$\mathbb{R}$. We formulate the computation in terms of the numbers
$t_{n}=(q^{(n+1)/2}-1)_{2}$ introduced at the end of
\S \ref{section:preliminaryresults}, where we recall from Lemma
\ref{lemma: t and w} that they are related to the numbers $w_{m}$ of Theorem
\ref{theorem2} by the formulae:%
\[
t_{8k+3}=w_{4k+2}\qquad\text{and}\qquad t_{8k+7}=w_{4k+4}\text{.}%
\]

\begin{theorem}
\label{theorem5} Up to finite groups of odd order, the groups ${}%
_{\varepsilon}KQ_{n}(\overline{R}_{F})$ are given in the following table.
(Recall that $\delta_{n0}$ denotes the Kronecker symbol.) \vspace{-0.1in}
\begin{table}[tbh]
\begin{center}%
\begin{tabular}
[c]{p{0.4in}|p{1.5in}|p{1.5in}|}\hline
$n\geq0$ & ${}_{-1}KQ_{n}(\overline{R}_{F})$ & ${}_{1}KQ_{n}(\overline{R}%
_{F})$\\\hline
$8k$ & $\delta_{n0}\mathbb{Z}$ & $\delta_{n0}\mathbb{Z}\oplus\mathbb{Z}%
\oplus\mathbb{Z}/2$\\
$8k+1$ & $0$ & $(\mathbb{Z}/2)^{3}$\\
$8k+2$ & $\mathbb{Z}$ & $(\mathbb{Z}/2)^{2}$\\
$8k+3$ & $\mathbb{Z}/2t_{8k+3}$ & $\mathbb{Z}/t_{8k+3}$\\
$8k+4$ & $\mathbb{Z}/2$ & $\mathbb{Z}$\\
$8k+5$ & $\mathbb{Z}/2$ & $0$\\
$8k+6$ & $\mathbb{Z}$ & $0$\\
$8k+7$ & $\mathbb{Z}/t_{8k+7}$ & $\mathbb{Z}/t_{8k+7}$\\\hline
\end{tabular}
\end{center}
\end{table}
\end{theorem}

\noindent\textbf{Proof.} Throughout the proof, we exploit Friedlander's
computation of ${}_{\varepsilon}KQ_{n}(\mathbb{F}_{q})$ given in \cite[Theorem
1.7]{Friedlander}, and work modulo odd torsion.\newline\noindent\textit{First
case}: $\varepsilon=1$.\ Applying Lemma \ref{_1KQ(R) splits} of Appendix B to
the homotopy fiber sequence (\ref{_epsilonKQ(R^{-}_F) fibration}) gives for
each $n$ a split short exact sequence%
\[
0\rightarrow{}_{1}KQ_{n}(\overline{R}_{F})\longrightarrow K_{n}(\mathbb{R}%
)\oplus K_{n}(\mathbb{R})\oplus{}_{1}KQ_{n}(\mathbb{F}_{q})\longrightarrow
K_{n}(\mathbb{R})\rightarrow0\text{.}%
\]
From this splitting we deduce an isomorphism
\[
{}_{1}KQ_{n}(\overline{R}_{F})\cong K_{n}(\mathbb{R})\oplus{}_{1}%
KQ_{n}(\mathbb{F}_{q}).
\]
\noindent\textit{Second case}: $\varepsilon=-1$.\ For the computation of
${}_{-1}\mathcal{KQ}(\overline{R}_{F})$ we need to make several case
distinctions arising from the $2$-completed homotopy cartesian square:%
\[%
\begin{array}
[c]{ccc}%
_{-1}\mathcal{KQ}(\overline{R}_{F})_{\#}^{c} & \longrightarrow &
_{-1}\mathcal{KQ}(\mathbb{R})_{\#}^{c}\simeq\mathcal{K}(\mathbb{C})_{\#}%
^{c}\simeq(\Omega^{4}(\mathcal{K}(\mathbb{C})))_{\#}^{c}\\
\downarrow &  & \downarrow\\
_{-1}\mathcal{KQ}(\mathbb{F}_{q})_{\#}^{c} & \longrightarrow & _{-1}%
\mathcal{KQ}(\mathbb{C})_{\#}^{c}\simeq\mathcal{K}(\mathbb{H})_{\#}^{c}%
\simeq(\Omega^{4}(\mathcal{K}(\mathbb{R})))_{\#}^{c}%
\end{array}
\]
Note that the vertical homotopy fiber is $(\Omega^{6}\mathcal{K}%
(\mathbb{R}))_{\#}^{c}$ by Lemma \ref{Connective.Spectra} since the map
$K_{4}(\mathbb{C)}\longrightarrow K_{4}(\mathbb{R)}$ is onto. The horizontal
homotopy fiber is $(\Omega^{5}\mathcal{K}(\mathbb{R}))_{\#}^{c}$ according to
\cite[Theorem 1.7]{Friedlander}. In particular, there is the \textquotedblleft
vertical\textquotedblright\ exact sequence
\begin{equation}
\cdots\rightarrow{}_{-1}KQ_{n+1}(\mathbb{F}_{q})\rightarrow K_{n+6}%
(\mathbb{R})\rightarrow{}_{-1}KQ_{n}(\overline{R}_{F})\rightarrow{}_{-1}%
KQ_{n}(\mathbb{F}_{q})\rightarrow K_{n+5}(\mathbb{R})\rightarrow\cdots,
\label{equation:vertical}%
\end{equation}
and the \textquotedblleft horizontal\textquotedblright\ exact sequence
\begin{equation}
\cdots\rightarrow K_{n+5}(\mathbb{C})\rightarrow K_{n+5}(\mathbb{R}%
)\rightarrow{}_{-1}KQ_{n}(\overline{R}_{F})\rightarrow K_{n+4}(\mathbb{C}%
)\rightarrow K_{n+4}(\mathbb{R})\rightarrow\cdots. \label{equation:horizontal}%
\end{equation}

If $n\equiv0,1\;(\mathrm{mod}\ 8)$ is nonzero, then (\ref{equation:vertical})
implies ${}_{-1}{KQ}_{n}(\overline{R}_{F})=0$. Likewise, when $n\equiv
2\;(\mathrm{mod}\ 8)$, (\ref{equation:vertical}) shows that ${}_{-1}{KQ}%
_{n}(\overline{R}_{F})\cong K_{8}(\mathbb{R})\cong\mathbb{Z}$.

For $n\equiv4\;(\mathrm{mod}\ 8)$, we use the segment
\[
\xymatrix{ K_{n+5}(\mathbb{C})\ar@{}[r]|-{\hbox{\large{$\rightarrow$}}}&
K_{n+5}(\mathbb{R})
\ar@{}[r]|-{\hbox{\large{$\rightarrow$}}}&{}_{-1}{KQ}_n(\overline{R}_{F})\ar@{}[r]|-{\hbox{\large{$\rightarrow$}}}&
K_{n+4}(\mathbb{C}) }
\]
of (\ref{equation:horizontal}). By analyzing (\ref{equation:vertical}), it
follows that ${}_{-1}{KQ}_{n}(\overline{R}_{F})$ is finite. Thus ${}_{-1}%
{KQ}_{n}(\overline{R}_{F})$ has order $2$.

For $n\equiv5\;(\mathrm{mod}\ 8)$, (\ref{equation:horizontal}) implies that
${}_{-1}{KQ}_{n}(\overline{R}_{F})$ is cyclic, whence by
(\ref{equation:vertical}) there is an isomorphism ${}_{-1}{KQ}_{n}%
(\overline{R}_{F})\cong\mathbb{Z}/2$.

For $n\equiv3\;(\mathrm{mod}\ 8)$, we use the exact sequence
(\ref{equation:vertical}) from ${}_{-1}{KQ}_{n+2}(\overline{R}_{F})$ to
$K_{n+5}(\mathbb{R})$. In view of the two previous results, this takes the
form%
\[
\mathbb{Z}/2\rightarrow\mathbb{Z}/2\oplus\mathbb{Z}/2\rightarrow
\mathbb{Z}/2\rightarrow\mathbb{Z}/2\rightarrow\mathbb{Z}/2\rightarrow
\mathbb{Z}/2\rightarrow{}_{-1}{KQ}_{n}(\overline{R}_{F})\rightarrow
\mathbb{Z}/t_{n}\rightarrow\mathbb{Z}\text{.}%
\]
Chasing this sequence from the left reveals that ${}_{-1}{KQ}_{n}(\overline
{R}_{F})$ is a finite group of order $2t_{n}$. Meanwhile, the exact sequence
(\ref{equation:horizontal}) obliges this group to be cyclic.

For $n\equiv6\;(\mathrm{mod}\ 8)$, (\ref{equation:vertical}) implies that
${}_{-1}{KQ}_{n}(\overline{R}_{F})\cong\mathbb{Z}$ or $\mathbb{Z}%
\oplus\mathbb{Z}/2$. On the other hand, the exact sequence
(\ref{equation:horizontal}) shows that ${}_{-1}{KQ}_{n}(\overline{R}_{F})$ is
a subgroup of $\mathbb{Z}$.

Finally, if $n\equiv7\;(\mathrm{mod}\ 8)$ then (\ref{equation:vertical})
produces an isomorphism between ${}_{-1}{KQ}_{n}(\overline{R}_{F})$ and
${}_{-1}{KQ}_{n}(\mathbb{F}_{q})\cong\mathbb{Z}/t_{n}$. \hfill$\Box\medskip$

We can now prove Theorem \ref{converse to theorem 1}.

\begin{theorem}
For every totally real number field $F$, the following are equivalent.

\begin{enumerate}
\item[(i)] $F$ is $2$-regular.

\item[(ii)] The square (\ref{hermitian K Bokstedt square}) is homotopy
cartesian for $F$ when $\varepsilon=1$.

\item[(iii)] The square (\ref{algebraic K Bokstedt square}) is homotopy
cartesian for $F$.
\end{enumerate}
\end{theorem}

\noindent\textbf{Proof.} Theorem \ref{theorem1} shows
that (i) implies (ii).

In the other direction, the preceding proof shows that (ii) leads to the
second column of the table of Theorem \ref{theorem5}. Since $W(R_{F})$ injects
into $W(F)$ \cite[IV (3.3)]{MH}, which has no odd-order torsion \cite[III
(3.10)]{MH}, we may work modulo odd torsion. From the epimorphisms ($i=0,1$)
${}_{1}{KQ}_{i}(\mathbb{R})\twoheadrightarrow{}_{1}{KQ}_{i}(\mathbb{C})$, we
obtain, from the Mayer-Vietoris sequence following from Theorem \ref{theorem1}%
, a short exact sequence
\[
0\rightarrow{}_{1}{KQ}_{0}(R_{F})\longrightarrow{}_{1}{KQ}_{0}(\mathbb{F}%
_{q})\oplus\bigoplus^{r}{}_{1}{KQ}_{0}(\mathbb{R})\longrightarrow\bigoplus
^{r}{}_{1}{KQ}_{0}(\mathbb{C})\rightarrow0\text{,}%
\]
which splits because the final group is free abelian. It follows that%
\[
{}_{1}{KQ}_{0}(R_{F})\cong\mathbb{Z}\oplus\mathbb{Z}^{r}\oplus\mathbb{Z}%
/2\text{,}%
\]
such that its $\mathbb{Z}/2$ summand maps nontrivially in the commuting square%
\[%
\begin{array}
[c]{ccccc}%
\mathbb{Z}\oplus\mathbb{Z}^{r}\oplus\mathbb{Z}/2\cong & {}_{1}{KQ}_{0}%
(R_{F}) & \longrightarrow & {}_{1}{KQ}_{0}(\mathbb{F}_{q}) & \cong%
\mathbb{Z}\oplus\mathbb{Z}/2\\
& \downarrow &  & \downarrow & \\
& {}_{1}W_{0}(R_{F}) & \longrightarrow & {}_{1}W_{0}(\mathbb{F}_{q}) &
\cong\mathbb{Z}/2\oplus\mathbb{Z}/2
\end{array}
\]
whose vertical maps are, by definition, surjective. Thus, in the exact
sequence%
\[
K_{0}(R_{F})\overset{H}{\longrightarrow}{}_{1}{KQ}_{0}(R_{F})\longrightarrow
{}_{1}W_{0}(R_{F})\rightarrow0
\]
the hyperbolic homomorphism $H$ must map the finite summand \textrm{Pic}%
$(R_{F})$ of $K_{0}(R_{F})$ trivially, and so have its cokernel $W(R_{F})$
isomorphic to $\mathbb{Z}^{r}\oplus\mathbb{Z}/2$. Then by Proposition
\ref{2+-regular characterization}(1), (7), it follows that $F$ is
$2$-regular.

Similarly, the $K$-theoretic theorem of \cite{HO} and \cite{Mitchell} in case
(i) asserts that (i) implies (iii).

Again, the computation $K_{2}(R_{F})\{2\}\cong(\mathbb{Z}/2)^{r}$ follows from
(iii). Here, Proposition \ref{2+-regular characterization}(1), (2) yield that
$F$ is $2$- regular. \hfill$\Box$


\section{Proof of Theorem \ref{theorem4}\label{V-computation}}

As for Theorem \ref{theorem2}, the splitting results proven in Lemmas
\ref{_1V(R_F) splits} and \ref{_-1V(R_F) splits} show that in order to prove
Theorem \ref{theorem4}, it suffices to compute the groups ${}_{\varepsilon
}V_{n}(\overline{R}_{F})$ introduced in the same section.

\begin{theorem}
\label{theorem6} Up to finite groups of odd order, the groups
\[
{}_{\varepsilon}V_{n}(\overline{R}_{F}):=\pi_{n}({}_{\varepsilon}%
\mathcal{V}(\overline{R}_{F}))
\]
are as follows. \begin{table}[tbh]
\begin{center}%
\begin{tabular}
[c]{p{0.4in}|p{1.0in}|p{1.2in}|}\hline
$n\geq0$ & ${}_{-1}V_{n}(\overline{R}_{F})$ & ${}_{1}V_{n}(\overline{R}_{F}%
)$\\\hline
$8k$ & $\mathbb{Z}\oplus\mathbb{Z}/2$ & $\mathbb{Z}^{2}$\\
$8k+1$ & $0$ & $(\mathbb{Z}/2)^{2}$\\
$8k+2$ & $\mathbb{Z}$ & $(\mathbb{Z}/2)^{2}$\\
$8k+3$ & $0$ & $0$\\
$8k+4$ & $\mathbb{Z}$ & $\mathbb{Z}^{2}$\\
$8k+5$ & $\mathbb{Z}/2$ & $0$\\
$8k+6$ & $\mathbb{Z}\oplus\mathbb{Z}/2$ & $0$\\
$8k+7$ & $\mathbb{Z}/2$ & $0$\\\hline
\end{tabular}
\end{center}
\end{table}

More precisely, for $\varepsilon=1$, the spectrum ${}_{1}\mathcal{V}%
(\overline{R}_{F})_{\#}^{c}$ (resp.~ ${}_{1}\mathcal{V}(R_{F})_{\#}^{c}$) has
the homotopy type of $2$ copies (resp.~ $2r$ copies) of the spectrum
$\mathcal{K}(\mathbb{R})_{\#}^{c}$.
\end{theorem}

\noindent\textbf{Proof.} Throughout the proof, we again work modulo odd
torsion.\newline\textit{First case}: $\varepsilon=1$.\ In the homotopy
cartesian square%
\[%
\begin{array}
[c]{ccc}%
{}_{1}\mathcal{V}(\overline{R}_{F})_{\#}^{c} & \longrightarrow & {}%
_{1}\mathcal{V}(\mathbb{R})_{\#}^{c}\simeq\mathcal{K}(\mathbb{R})_{\#}^{c}\\
\downarrow &  & \downarrow^{{}_{1}\sigma}\\
{}_{1}\mathcal{V}(\mathbb{F}_{q})_{\#}^{c} & \overset{\chi}{\longrightarrow} &
{}_{1}\mathcal{V}(\mathbb{C})_{\#}^{c}\simeq\Omega^{-1}\mathcal{K}%
(\mathbb{R})_{\#}^{c}%
\end{array}
\]
as noted in Lemma \ref{_1VR -> C} of Appendix B, the map ${}_{1}\sigma$ is
nullhomotopic. If we replace $\chi$ by a Serre fibration, then ${}%
_{1}\mathcal{V}(\overline{R}_{F})_{\#}^{c}$ has the homotopy type of its
pullback over the nullhomotopic map ${}_{1}\sigma$. This means that ${}%
_{1}\mathcal{V}(\overline{R}_{F})_{\#}^{c}$ has the homotopy type of a fiber
homotopy trivial fibration with base ${}_{1}\mathcal{V}(\mathbb{R})_{\#}%
^{c}=\mathcal{K}(\mathbb{R})_{\#}^{c}$, and fiber of the homotopy type of the
fiber of the lower horizontal map, which has been\textsf{ }identified in
\cite[Corollary 1.6]{Friedlander} as $\mathcal{K}(\mathbb{R})_{\#}^{c}$. By
Lemmas \ref{_1V(R_F) splits} and \ref{_-1V(R_F) splits}, a similar argument
holds for $_{1}\mathcal{V}(R_{F})_{\#}^{c}$. Hence, the spectrum ${}%
_{1}\mathcal{V(}\overline{R}_{F})_{\#}^{c}$ (resp${}$ $_{1}\mathcal{V}%
(R_{F})_{\#}^{c}$) has the homotopy type of $\vee^{2}\mathcal{K}%
(\mathbb{R})_{\#}^{c}$ (resp.~ $\vee^{2r}\mathcal{K}(\mathbb{R})_{\#}^{c}$).

\textit{Second case}: $\varepsilon=-1$.\ Combining the results of Friedlander
\cite[Corollary 1.6]{Friedlander} and Quillen \cite{Quillen} via the forgetful
map relating algebraic and hermitian $K$-theory gives a homotopy commutative
diagram%
\[%
\begin{array}
[c]{ccccc}%
{}_{-1}\mathcal{V}(\mathbb{F}_{q})_{\#}^{c} & \longrightarrow & {}%
_{-1}\mathcal{V}(\mathbb{C})_{\#}^{c} & \longrightarrow & {}_{-1}%
\mathcal{V}(\mathbb{C})_{\#}^{c}\\
\downarrow &  & \downarrow &  & \downarrow\\
{}_{-1}\mathcal{KQ}(\mathbb{F}_{q})_{\#}^{c} & \longrightarrow & {}%
_{-1}\mathcal{KQ}(\mathbb{C})_{\#}^{c} & \overset{\psi^{q}-1}{\longrightarrow}
& {}_{-1}\mathcal{KQ}(\mathbb{C})_{\#}^{c}\\
\downarrow &  & \downarrow &  & \downarrow\\
\mathcal{K}(\mathbb{F}_{q})_{\#}^{c} & \longrightarrow & \mathcal{K}%
(\mathbb{C})_{\#}^{c} & \overset{\psi^{q}-1}{\longrightarrow} & \mathcal{K}%
(\mathbb{C})_{\#}^{c}%
\end{array}
\]
in which all colums and both lower rows are homotopy fibrations/cofibrations.
It follows from \cite[Lemma 2.1]{CMN} once more that the top row is also a
fibration/cofibration, with the cofiber map induced from $\psi^{q}-1$. It
follows that the defining homotopy cartesian square for ${}_{-1}%
\mathcal{V}(\overline{R}_{F})_{\#}^{c}$ from Section
\ref{section:splittingresults}, by means of substitutions according to Lemma
\ref{_-1VRC} in Appendix B, gives rise to a homotopy commutative diagram where
the horizontal maps are homotopy fibrations (note that $\pi_{0}(\Omega
^{3}\mathcal{K}(\mathbb{R}))=0$):%
\begin{equation}%
\begin{array}
[c]{ccccc}%
{}_{-1}\mathcal{V}(\overline{R}_{F})_{\#}^{c} & \longrightarrow & {}%
_{-1}\mathcal{V}(\mathbb{R})_{\#}^{c}\simeq(\Omega^{2}\mathcal{K}%
(\mathbb{R}))_{\#}^{c} & \overset{\tau}{\longrightarrow} & (\Omega
^{3}\mathcal{K}(\mathbb{R}))_{\#}^{c}\\
\downarrow &  & \downarrow^{{}_{-1}\sigma} &  & \downarrow^{\mathrm{id}}\\
{}_{-1}\mathcal{V}(\mathbb{F}_{q})_{\#}^{c} & \longrightarrow & {}%
_{-1}\mathcal{V}(\mathbb{C})_{\#}^{c}\simeq(\Omega^{3}\mathcal{K}%
(\mathbb{R}))_{\#}^{c} & \overset{\psi^{q}-1}{\longrightarrow} & (\Omega
^{3}\mathcal{K}(\mathbb{R}))_{\#}^{c}%
\end{array}
\label{-1VDiagram}%
\end{equation}

From the proof of Lemma \ref{_-1V(R) -> _-1V(C)} in Appendix B again, the map
${}_{-1}\sigma$ corresponds to the cup-product with the generator of
$K_{1}(\mathbb{R})$. Now, the image of $_{-1}\sigma$ is a torsion element in
$K_{\ast}(\mathbb{R})$; and it is easy to see by a direct checking that, for
odd $q$, any torsion element of $K_{\ast}(\mathbb{R})$ is killed by $\psi
^{q}-1$. Therefore, the homotopy exact sequence for the upper horizontal maps
includes the exact sequence%
\[
\xymatrix{ K_{n+4}(\mathbb{R}) \ar@{}[r]|-{\hbox{\large{$\rightarrow$}}}&
{}_{-1}V_{n}(\overline{R}_{F}) \ar@{}[r]|-{\hbox{\large{$\rightarrow$}}}&
K_{n+2}(\mathbb{R}) \ar@{}[r]|-{\hbox{\large{$\rightarrow$}}}& K_{n+3}(\mathbb{R}), }
\]
in which the last map is trivial because it is given by the cup-product with
the generator of $K_{1}(\mathbb{R})$ composed with $\psi^{q}-1$. The groups
${}_{-1}V_{n}(\overline{R}_{F})$ are therefore included in the short exact
sequences%
\[
0\longrightarrow K_{n+4}(\mathbb{R})\longrightarrow{}_{-1}V_{n}(\overline
{R}_{F})\longrightarrow K_{n+2}(\mathbb{R})\longrightarrow0
\]
which determine them, except for $n\equiv0\;(\mathrm{mod}\ 8)$, where we can
say only that ${}_{-1}V_{n}(\overline{R}_{F})=\mathbb{Z}$ or $\mathbb{Z}%
/2\oplus\mathbb{Z}$. One way to resolve this ambiguity is to write the exact
sequence%
\[
0={}_{-1}KQ_{n+1}(\overline{R}_{F})\longrightarrow K_{n+1}(\overline{R}%
_{F})\longrightarrow{}_{-1}V_{n}(\overline{R}_{F})\longrightarrow{}_{-1}%
KQ_{n}(\overline{R}_{F})=0
\]
which implies that $K_{n+1}(\overline{R}_{F})\cong{}_{-1}V_{n}(\overline
{R}_{F})$.

In general, the computation of the groups $K_{n}(\overline{R}_{F})$ for $n>0$
follows from the analog of Diagram ($\ref{RFoverbar}$) for $\mathcal{K}%
(\overline{R}_{F})$. They are the following for $n\equiv k\;(\mathrm{mod}%
\ 8)$, starting from $k=0$:%
\[
0,\ \mathbb{Z}/2\oplus\mathbb{Z},\ \mathbb{Z}/2,\ \mathbb{Z}/2w_{4k+2}%
,\ 0,\ \mathbb{Z},\ 0,\ \mathbb{Z}/w_{4k+1}\text{.}%
\]

This computation is straightforward, except for $n\equiv1\;(\mathrm{mod}\ 8)$,
where we have to use two exact sequences extracted from the analog of the
previous square for the spectrum \ $\mathcal{K}(\overline{R}_{F})$. The first
one%
\[
0\longrightarrow K_{n+1}(\mathbb{C})\longrightarrow K_{n}(\overline{R}%
_{F})\longrightarrow K_{n}(\mathbb{R})\longrightarrow K_{n}(\mathbb{C})=0
\]
shows as expected that $K_{n}(\overline{R}_{F})=$ $\mathbb{Z}$ or
$\mathbb{Z}/2\oplus\mathbb{Z}$. In the second one, we write the Mayer-Vietoris
exact sequence associated to the same previous square:%
\[
0\longrightarrow K_{n+1}(\mathbb{C})=\mathbb{Z}\longrightarrow K_{n}%
(\overline{R}_{F})\longrightarrow K_{n}(\mathbb{R})\oplus K_{n}(\mathbb{F}%
_{q})\longrightarrow K_{n}(\mathbb{C})=0\text{.}%
\]
It shows that $K_{n}(\overline{R}_{F})=\pi_{n}(\mathcal{K}(\overline{R}_{F}))$
cannot be isomorphic to $\mathbb{Z}$, since $K_{n}(\mathbb{R})\oplus
K_{n}(\mathbb{F}_{q})$ is a direct sum of two nontrivial cyclic groups. The
computation of the groups ${}_{-1}V_{n}(\overline{R}_{F})$ is therefore
accomplished for all values of $n$.$\hfill\Box\smallskip$

\begin{remark}
On the level of spectra the composition $\tau$ in (\ref{-1VDiagram}),
\[
{}_{-1}\mathcal{V}(\mathbb{R})_{\#}^{c}\sim(\Omega^{2}\mathcal{K}%
(\mathbb{R}))_{\#}^{c}\overset{\sigma}{\longrightarrow}(\Omega^{3}%
\mathcal{K}(\mathbb{R}))_{\#}^{c}\overset{\Omega^{3}(\psi^{q}-1)}%
{\longrightarrow}(\Omega^{3}\mathcal{K}(\mathbb{R}))_{\#}^{c}%
\]
where $\sigma$ is induced by the cup-product with the generator of
$K_{1}(\mathcal{\mathbb{R}}),$ is NOT nullhomotopic. This fact is proved in
Appendix D.
\end{remark}

We can now use Theorems \ref{theorem2} and \ref{theorem4} to determine the
composition%
\[
{}_{\varepsilon}{KQ}_{n}(R_{F})\overset{F}{\longrightarrow}{K}_{n}%
(R_{F})\overset{H}{\longrightarrow}{}_{\varepsilon}{KQ}_{n}(R_{F})
\]
of the homomorphisms induced by the forgetful and hyperbolic functors. From
their respective induced homotopy fiber sequences ${}_{\varepsilon}%
\mathcal{V}(R_{F})\longrightarrow{}_{\varepsilon}\mathcal{KQ}(R_{F}%
)\longrightarrow\mathcal{K}(R_{F})$, ${}_{\varepsilon}\mathcal{U}%
(R_{F})\longrightarrow\mathcal{K}(R_{F})\longrightarrow{}_{\varepsilon
}\mathcal{KQ}(R_{F})$, and the natural homotopy equivalence ${}_{\varepsilon
}\mathcal{V}(R_{F})\simeq\Omega{}_{-\varepsilon}\mathcal{U}(R_{F})$ of
\cite{K:AnnM112fun}, we have the exact sequences
\[
\cdots\rightarrow{}_{\varepsilon}V_{n}(R_{F})\longrightarrow{}_{\varepsilon
}KQ_{n}(R_{F})\overset{F}{\longrightarrow}K_{n}(R_{F})\longrightarrow
{}_{\varepsilon}V_{n-1}(R_{F})\rightarrow\cdots,
\]
and%
\[
\cdots\rightarrow{}_{-\varepsilon}V_{n-1}(R_{F})\longrightarrow K_{n}%
(R_{F})\overset{H}{\longrightarrow}{}_{\varepsilon}KQ_{n}(R_{F}%
)\longrightarrow{}_{-\varepsilon}V_{n-2}(R_{F})\rightarrow\cdots\text{.}%
\]
Since all terms are now known (and many are zero), a routine computation gives
the following.

\begin{corollary}
\label{HF computation}For $n\geq1$, the endomorphism $HF$ of the group
${}_{\varepsilon}{KQ}_{n}(R_{F})$ modulo odd torsion

\begin{enumerate}
\item[(i)] is multiplication by $2$, when $n\equiv3\;(\mathrm{mod}\ 4)$
($\varepsilon=\pm1$),

\item[(ii)] has image of order $2$, when both $n\equiv1,2\;(\mathrm{mod}\ 8)$
and $\varepsilon=1$, and

\item[(iii)] is zero otherwise.$\hfill\Box\smallskip$
\end{enumerate}
\end{corollary}

A similar computation affords the corresponding result for the other
composition of $F$ and $H$. From \cite[pg.~230]{K:AnnM112hgo} we note that
this endomorphism of ${K}_{n}(R_{F})$ is the sum of the identity and the
involution induced by the duality functor. Since this involution is
independent of $\varepsilon$, we need consider only the simpler case
$\varepsilon=-1$.

\begin{corollary}
For $n\geq1$ and $\varepsilon=\pm1$, the endomorphism $FH$ of ${K}_{n}(R_{F})$
modulo odd torsion

\begin{enumerate}
\item[(i)] is multiplication by $2$, when $n\equiv3\;(\mathrm{mod}\ 4)$, and

\item[(ii)] is zero, otherwise.$\hfill\Box\smallskip$
\end{enumerate}
\end{corollary}

\begin{corollary}
The canonical involution on $K_{n}(R_{F})$ modulo odd torsion

\begin{enumerate}
\item[(i)] is the identity, for $n=0$ and $n\equiv3\;(\mathrm{mod}\ 4)$, and

\item[(ii)] is the opposite of the identity, otherwise.$\hfill\Box\smallskip$
\end{enumerate}
\end{corollary}

\begin{remark}
Concerning the odd torsion, in general the functors $F$ and $H$ induce
bijections between the symmetric parts of the $K$- and $KQ$-groups of a ring
$A$. Here the involution on the $K$-groups is induced by the duality functor.
Clearly the composition $FH$ is the multiplication by $2$ map on the symmetric
part. The same result holds for $HF$ on the symmetric part, while it is
trivial on the antisymmetric part. We note that it remains to compute the odd
torsion part of ${}_{\varepsilon}KQ_{n}(R_{F})$, even for $F=\mathbb{Q}$.
However, as we have seen more generally in Proposition \ref{Wodd}, the odd
torsion of the higher Witt groups and coWitt groups of $R_{F}$ is trivial.
\end{remark}

\section{Proof of Theorem \ref{theorem3}}

\label{section:proofoftheorem3}

In the terminology at the end of the Introduction, consider the naturally
induced map between $2$-completed connective spectra induced by the forgetful
functor and the homotopy fixed point functor for the $_{\varepsilon}%
\mathbb{Z}/2$-action:
\begin{equation}%
\begin{array}
[c]{ccc}%
{}_{\varepsilon}\mathcal{KQ}(R_{F})_{\#}^{c} & \rightarrow &
{\displaystyle \bigvee\limits^{r}}{}_{\varepsilon}\mathcal{KQ}(\mathbb{R}%
)_{\#}^{c}\\
\downarrow &  & \downarrow\\
{}_{\varepsilon}\mathcal{KQ}(\mathbb{F}_{q})_{\#}^{c} & \rightarrow &
{\displaystyle\bigvee\limits^{r}}{}_{\varepsilon}\mathcal{KQ}(\mathbb{C}%
)_{\#}^{c}%
\end{array}
\quad\rightarrow\quad%
\begin{array}
[c]{ccc}%
\mathcal{K}(R_{F})_{\#}{}^{h({}_{\varepsilon}\mathbb{Z}/2)} & \rightarrow &
{\displaystyle\bigvee\limits^{r}}\mathcal{K}(\mathbb{R})_{\#}^{c}{}%
^{h({}_{\varepsilon}\mathbb{Z}/2)}\\
\downarrow &  & \downarrow\\
\mathcal{K}(\mathbb{F}_{q})_{\#}{}^{h({}_{\varepsilon}\mathbb{Z}/2)} &
\rightarrow & {\displaystyle\bigvee\limits^{r}}\mathcal{K}(\mathbb{C}%
)_{\#}^{c}{}^{h({}_{\varepsilon}\mathbb{Z}/2)}%
\end{array}
\label{homotopycartesiandiagrams1}%
\end{equation}
As noted in the beginning of Section \ref{section:proofoftheorem1}, the
spectrum maps in the B{\"{o}}kstedt square are $_{\varepsilon}\mathbb{Z}%
/2$-equivariant, being induced by ring maps. Theorem \ref{theorem1} and the
main results in \cite{HO}, \cite{Mitchell} (\textsl{cf.}~Appendix
\ref{section:K-theorybackground} for more details) show that both the
hermitian and the algebraic $K$-theory squares are homotopy cartesian squares
(since the homotopy fixed point functor is a homotopy functor).

By \cite[Lemmas 7.3-7.5]{BK} the map
\[
{}_{\varepsilon}\mathcal{KQ}(A)_{\#}^{c}\rightarrow(\mathcal{K(A)}_{\#}%
{}^{h({}_{\varepsilon}\mathbb{Z}/2)})^{c}%
\]
in (\ref{homotopycartesiandiagrams1}) is a homotopy equivalence for
$A=\mathbb{F}_{q},\mathbb{R},\mathbb{C}$. It is worth mentioning that the most
delicate case is when $A=\mathbb{R},$ where the machinery of Fredholm
operators in an infinite dimensional real Hilbert space is used. It follows
that the induced map of homotopy pullbacks is also a homotopy
equivalence.\hfill$\Box\smallskip$

\appendix

\section{$K$-theory background}

\label{section:K-theorybackground}

In this appendix we deduce the homotopy cartesian square of $K$-theory spectra
(\ref{algebraic K Bokstedt square}) using the space level results given in
\cite{HO}. The examples of \'{e}tale $K$-theory spectra of real number fields
at the prime $2$ in \cite[\S 5]{Mitchell} provide an alternate proof on
account of the solution of the Quillen-Lichtenbaum conjecture in
\cite{Ostvar}. Throughout we retain the assumptions and notations employed in
the main body of the text. Recall from the Introduction that $q$ is a prime
number. \vspace{0.1in}

Recall the construction of the square (\ref{algebraic K Bokstedt square}) from
the beginning of Section \ref{section:proofoftheorem1}: one starts out by
choosing an embedding of the field of $q$-adic numbers $\mathbb{Q}_{q}$ into
the complex numbers $\mathbb{C}$ such that the induced composite map
\[
\mathcal{K}(\mathbb{Z}_{q})_{\#}\rightarrow\mathcal{K}(\mathbb{Q}_{q}%
)_{\#}\rightarrow\mathcal{K}(\mathbb{C})_{\#}^{c}
\]
agrees with Quillen's Brauer lift $\mathcal{K}(\mathbb{F}_{q})_{\#}%
\rightarrow\mathcal{K}(\mathbb{C})_{\#}^{c}$ from \cite{Quillen} under the
rigidity equivalence between $\mathcal{K}(\mathbb{Z}_{q})_{\#}$ and
$\mathcal{K}(\mathbb{F}_{q})_{\#}$ \cite{Gabber}. The ring maps relating
$R_{F}$ to $\mathbb{F}_{q}$, $\mathbb{R}$ and $\mathbb{C}$ induce the
commuting B{\"{o}}kstedt square (\ref{algebraic K Bokstedt square}) via
Suslin's identifications of the $2$-completed algebraic $K$-theory spectra of
the real numbers with $\mathcal{K}(\mathbb{R})_{\#}^{c}$ and likewise for the
complex numbers and $\mathcal{K}(\mathbb{C})_{\#}^{c}$ \cite{Suslin :local
fields}.

\begin{theorem}
The B{\"{o}}kstedt square (\ref{algebraic K Bokstedt square}) is a commuting
homotopy cartesian square of $2$-completed spectra:
\end{theorem}

\noindent\textbf{Proof.} Let $JK(q)$ denote the fiber of the composite map
\[
BO_{\#}\overset{c}{\longrightarrow}BU_{\#}\overset{\psi^{q}-1}{\longrightarrow
}BU_{\#}%
\]
where $c$ denotes the complexification map and $\psi^{q}$ the $q$\thinspace th
Adams operation on the $2$-completion of the classifying space $BU$. As usual,
$U$ and $O$ are the stable unitary and orthogonal groups. When $q\equiv
\pm3\;(\mathrm{mod}\ 8)$, $JK(q)$ is a space level model for $\mathcal{K}%
(R_{\mathbb{Q}})_{\#}$ \cite{HO}. By the main result in \cite{Quillen}, the
fiber of $\psi^{q}-1$ identifies with the $2$-completed algebraic $K$-theory
space of $\mathbb{F}_{q}$. Moreover, the product decomposition
\[
JK(q)\times\prod^{r-1}U_{\#}/O_{\#}%
\]
of the $2$-completed algebraic $K$-theory space of $R_{F}$ established in
\cite[Theorem 1.1]{HO} shows that the space level analogue of
(\ref{algebraic K Bokstedt square}) is homotopy cartesian. On the other hand,
the Quillen-Lichtenbaum conjecture for totally real number fields
\cite{Ostvar} implies that $\mathcal{K}(R_{F})_{\#}$ is homotopy equivalent to
the connective cover of its $K(1)$-localization $L_{K(1)}\mathcal{K}%
(R_{F})_{\#}$, and likewise for $\mathbb{F}_{q}$, $\mathbb{R}$ and
$\mathbb{C}$. Here $K(1)$ is the first Morava $K$-theory spectrum at the prime
$2$. In order to conclude we incorporate \cite{Bousfield}, which reduces
questions about $K(1)$-local spectra to space level questions. That is,
applying Bousfield's homotopy functor $T$ from spaces to spectra yields the
desired conclusion since by \textsl{loc.~cit.}~$L_{K(1)}\mathcal{K}%
(R_{F})_{\#}$ identifies with $T\Omega^{\infty}\mathcal{K}(R_{F})_{\#}$.
\hfill$\Box$\smallskip

We refer the reader to \cite{MitchellSurvey} for an extensive background on
the stable homotopy-theoretic interpretation of the Quillen-Lichtenbaum conjecture.

\begin{remark}
As we already mentioned in the Introduction (see Remark \ref{localization}),
we also have a \textquotedblleft B{\"{o}}kstedt square\textquotedblright\ if
we decide to consider $2$-localizations instead of $2$-completions, provided
we follow the convention of our paper that the fields $\mathbb{R}$ and
$\mathbb{C}$ are considered with their usual topology.
\end{remark}

\section{Homology module maps}

Most theories in this paper are modules over the graded ring ${}_{\varepsilon
}{KQ}_{\ast}(R_{F})$ in the case $F=\mathbb{Q}$, in other words,
${}_{\varepsilon}{KQ}_{\ast}(\mathbb{Z}[1/2])$. The framework for such
considerations is laid out in \cite[\S 3]{K:asterisque149}, using the
description of algebraic $K$-theory in terms of flat \textquotedblleft
virtual\textquotedblright\ bundles.

In the topological case when $A=\mathbb{R}$, $\mathbb{C}$, the module
structures on $\mathcal{K}(A)$, ${}_{\varepsilon}\mathcal{KQ}(A)$,
${}_{\varepsilon}\mathcal{U}(A)$ and ${}_{\varepsilon}\mathcal{V}(A)$ are much
simpler to define. For clarity, we discuss the examples of ${}_{\varepsilon
}\mathcal{V}(\mathbb{R})$ and ${}_{\varepsilon}\mathcal{V}(\mathbb{C})$,
leading to a determination of the homotopy fiber of the map
\[
{}_{\varepsilon}\mathcal{V}(\mathbb{R})\longrightarrow{}_{\varepsilon
}\mathcal{V}(\mathbb{C})
\]
for both cases $\varepsilon=\pm1$.

We start with a geometric viewpoint: the cohomology theory associated to the
spectrum ${}_{1}\mathcal{KQ}(\mathbb{R})$ is constructed as the $K$-theory of
real vector bundles equipped with nondegenerate quadratic forms. As shown in
\cite[Exercise 9.22]{K:ktheorybook}, such a vector bundle $E$ splits as a
Whitney sum
\[
E=E^{+}\oplus E^{-},
\]
where the quadratic form is positive on $E^{+}$ and negative on $E^{-}$. A
bundle version of Sylvester's theorem tells us that the isomorphism classes of
$E^{+}$ and $E^{-}$ are independent of the sum decomposition (see also the
remarks below). Hence, we have the following assertion.

\begin{lemma}
\label{_1KQ(R) splits}There are splittings
\[
{}_{1}\mathcal{KQ}(\mathbb{R})\simeq\mathcal{K}(\mathbb{R})\vee\mathcal{K}%
(\mathbb{R})\text{\quad and\quad}_{1}\mathcal{KQ}(\mathbb{R})^{c}%
\simeq\mathcal{K}(\mathbb{R})^{c}\vee\mathcal{K}(\mathbb{R})^{c}\text{.}%
\]
Moreover, the $K$-theory of the category of real vector bundles, on a compact
space $X$ and provided with a nondegenerate quadratic form, is canonically
isomorphic to the direct sum of two copies of the usual real $K$-theory of
$X$.$\hfill\Box\smallskip$
\end{lemma}

Another result of interest is also shown in \cite[Exercise 9.22]%
{K:ktheorybook}:

\begin{lemma}
\label{_1KQ(C)}There are homotopy equivalences of spectra%
\[
{}_{1}\mathcal{KQ}(\mathbb{C})\simeq\mathcal{K}(\mathbb{R})\text{\quad
and\quad}_{1}\mathcal{KQ}(\mathbb{C})^{c}\simeq\mathcal{K}(\mathbb{R}%
)^{c}\text{.}%
\]
Moreover, the $K$-theory of the category of complex vector bundles, on a
compact space $X$ and provided with a nondegenerate quadratic form, is
canonically isomorphic to the usual real $K$-theory of $X$.
\end{lemma}

\begin{remark}
A slightly different proof of these lemmas is to use the following classical
result: a Lie group has the homotopy type of its compact form \cite[pp.
218--219]{Helgason}. For instance, the Lie group $O(p,q)$ has the homotopy
type of $O(p)\times O(q)$. This implies that the homotopy theory of real
vector bundles provided with a quadratic form of type $(p,q)$ is equivalent to
the homotopy theory of couples of real vector bundles $(E^{+},E^{-})$, of
dimensions $p$ and $q$ respectively. A similar example of interest is the Lie
group $O(n,\mathbb{C})$ which has the homotopy type of the usual compact Lie
group $O(n)$. This implies that the homotopy theory of complex vector bundles
of dimension $n$ provided with a nondegenerate quadratic form is equivalent to
the homotopy theory of real vector bundles of rank $n$. These results of
course imply the previous lemmas. Moreover, since all theories involved are
$8$-periodic according to Bott, the homotopy equivalences on the level of the
$0$-space imply the homotopy equivalence of spectra.
\end{remark}

In the considerations that follow, we prefer to take the bundle viewpoint
which is easier to handle than its homotopy counterpart, especially for module
or ring structures which are simply given by the tensor product of vector
bundles in the appropriate categories.

We illustrate this philosophy by a concrete description of the spectrum
${}_{1}\mathcal{V}(\mathbb{R})$ which is the homotopy fiber of the forgetful
map ${}_{1}\mathcal{KQ}(\mathbb{R})\overset{F}{\rightarrow}\mathcal{K}%
(\mathbb{R})$. Strickly speaking, one should describe the full spectrum.
However, by classical Bott periodicity, it is enough to describe the $0$-part
of the spectrum. Since ${}_{1}\mathcal{KQ}(\mathbb{R})$ splits as
$\mathcal{K}(\mathbb{R})\vee\mathcal{K}(\mathbb{R})$, the map $F$ being
induced by the direct sum, the homotopy fiber should be $\mathcal{K}%
(\mathbb{R})$. We want to be more precise in terms of module structures and
consider the \textquotedblleft relative\textquotedblright\ cohomology theory
associated to this homotopy fiber ${}_{1}\mathcal{V}(\mathbb{R})$. It can be
described by a well-known scheme going back to Atiyah-Hirzebruch
\cite{Atiyah-Hirzebruch}, reproduced in \cite[pp. 59--63]{K:ktheorybook} (in a
slightly different context) and also in \cite[pg $269$]{K:AnnM112fun}. One
considers homotopy classes of triples $\tau=(E,F,\alpha)$, where $E$ and $F$
are real vector bundles equipped with nondegenerate quadratic forms, and
$\alpha$ is an isomorphism between the underlying real vector bundles. If $G$
is another real vector bundle with a nondegenerate quadratic form, then its
cup-product with $\tau$ is given as the triple
\[
(G\otimes E,\,G\otimes F,\,\mathrm{id}\otimes\alpha)\text{.}%
\]
This defines a ${}_{1}\mathcal{KQ}(\mathbb{R})$-module structure on ${}%
_{1}\mathcal{V}(\mathbb{R})$. By associating to every real vector bundle a
metric,\textsl{ i.e}.~a positive quadratic form, we obtain a well-defined map
up to homotopy $\mathcal{K}(\mathbb{R})\rightarrow{}_{1}\mathcal{KQ}%
(\mathbb{R})$, which is a right inverse to the forgetful map. Therefore, every
${}_{1}\mathcal{KQ}(\mathbb{R})$-module acquires a naturally induced
$\mathcal{K}(\mathbb{R})$-module structure.

\begin{lemma}
\label{_1VR}The spectrum ${}_{1}\mathcal{V}(\mathbb{R})$ is homotopy
equivalent to the real topological $K$-theory $\mathcal{K}(\mathbb{R})$ as a
$\mathcal{K}(\mathbb{R})$-module spectrum, and hence $_{1}\mathcal{V}%
(\mathbb{R})^{c}\simeq\mathcal{K}(\mathbb{R})^{c}$ as $\mathcal{K}%
(\mathbb{R})^{c}$-module spectra.
\end{lemma}

\noindent\textbf{Proof.} We can identify ${}_{1}\mathcal{V}(\mathbb{R})$ with
$\mathcal{K}(\mathbb{R})$ as modules over $\mathcal{K}(\mathbb{R})$ as
follows: if $E$ is a real vector bundle there is an associated triple
$(E_{+},E_{-},\mathrm{id})$ where $E_{+}$ is the bundle $E$ equipped with a
positive quadratic form, and likewise for $E_{-}$ but with a negative
quadratic form. This correspondence has an inverse defined by associating to a
triple $\tau=(E,F,\alpha)$ as above the formal difference $E_{+}-F_{+}$ of the
respective positive-form summands.\hfill$\Box$\smallskip

\begin{lemma}
\label{_1VC}The spectrum ${}_{1}\mathcal{V}(\mathbb{%
\mathbb{C}
})$ is homotopy equivalent to the spectrum $\Omega^{-1}(\mathcal{K}%
(\mathbb{R}))$, and therefore ${}_{1}\mathcal{V}(\mathbb{%
\mathbb{C}
})^{c}\simeq\Omega^{-1}(\mathcal{K}(\mathbb{R})^{c})$.
\end{lemma}

\noindent\textbf{Proof.} The theory ${}_{1}\mathcal{V}(\mathbb{C})$ is the
homotopy fiber of the map ${}_{1}\mathcal{KQ}(\mathbb{C})\overset
{F}{\longrightarrow}\mathcal{K}(\mathbb{C})$, which arises from triples
$(E_{1},E_{2},\alpha)$, where $E_{1}$ and $E_{2}$ are real vector bundles and
$\alpha$ is an isomorphism between their corresponding complexified vector
bundles. By a well known theorem of Bott (see for instance \cite[Section
III.5]{K:ktheorybook}) , this homotopy fiber may be identified with
$\Omega^{-1}(\mathcal{K}(\mathbb{R}))$. Since $K_{-1}(\mathbb{R})=0$, we also
have ${}_{1}\mathcal{V}(\mathbb{%
\mathbb{C}
})^{c}\simeq\Omega^{-1}(\mathcal{K}(\mathbb{R})^{c})$ according to Lemma
\ref{Connective.Spectra}.\hfill$\Box\smallskip$

\begin{lemma}
\label{_1VR -> C}The map ${}_{1}\mathcal{V}(\mathbb{R})\longrightarrow{}%
_{1}\mathcal{V}(\mathbb{C})$ is nullhomotopic and its homotopy fiber has the
homotopy type of $\mathcal{K}(\mathbb{R})\vee\mathcal{K}(\mathbb{R})$. In the
same way, the map ${}_{1}\mathcal{V}(\mathbb{R})^{c}\longrightarrow{}%
_{1}\mathcal{V}(\mathbb{C})^{c}$ is nullhomotopic and its homotopy fiber has
the homotopy type of $\mathcal{K}(\mathbb{R})^{c}\vee\mathcal{K}%
(\mathbb{R})^{c}$.
\end{lemma}

\noindent\textbf{Proof.} By the above considerations, the two theories also
have a $\mathcal{K}(\mathbb{R})$-module structure. Since ${}_{1}%
\mathcal{V}(\mathbb{R})$ is free of rank one as a $\mathcal{K}(\mathbb{R}%
)$-module, the map ${}_{1}\mathcal{V}(\mathbb{R})\longrightarrow{}%
_{1}\mathcal{V}(\mathbb{C})$ is determined up to homotopy equivalence by its
effect on the zeroth homotopy groups. This means that the map of associated
real topological $K$-theories%
\[
K_{\mathbb{R}}^{n}(X)\longrightarrow K_{\mathbb{R}}^{n+1}(X)
\]
is induced by the cup-product with an element of $K_{\mathbb{R}}%
^{1}(\mathrm{point})=K_{-1}(\mathbb{R})=0$. This shows that our first map is
nullhomotopic. Its fiber has the homotopy type of $\mathcal{K}(\mathbb{R}%
)\vee\mathcal{K}(\mathbb{R})$ since the fiber of this nullhomotopic fibration
has the homotopy type of the product of the total space (see (\ref{_1VR}))
with the loop space of the base (see (\ref{_1VC})). The same statements hold
for connective covers since, from the previous lemma for example, the map
${}_{1}V_{0}(%
\mathbb{R}
)\longrightarrow{}_{1}V_{0}(%
\mathbb{C}
)=0$ is an epimorphism.\hfill$\Box\smallskip$

The determination of the map ${}_{-1}\mathcal{V}(\mathbb{R})\rightarrow{}%
_{-1}\mathcal{V}(\mathbb{C})$ is more delicate but uses the same arguments. By
definition, the spectrum $_{-1}\mathcal{V}(\mathbb{R})$ is the fiber of
${}_{-1}\mathcal{KQ}(\mathbb{R})=\mathcal{K}(\mathbb{C})\rightarrow
\mathcal{K}(\mathbb{R})$ and so identifies with $\Omega^{2}(\mathcal{K}%
(\mathbb{R}))$, according to a classical result of Bott. For the same reasons,
$_{-1}\mathcal{V}(\mathbb{R})$, which is the fiber of%
\[
_{-1}\mathcal{KQ}(\mathbb{C})\simeq\Omega^{4}(\mathcal{K}(\mathbb{R}%
))\longrightarrow\mathcal{K}(\mathbb{C})\simeq\Omega^{4}(\mathcal{K}%
(\mathbb{C}))\text{,}%
\]
is homotopically equivalent to $\Omega^{3}(\mathcal{K}(\mathbb{R}))$.
Summarizing, we have proved the following lemma:

\begin{lemma}
\label{_-1VRC}The spectrum ${}_{-1}\mathcal{V}(\mathbb{R})$ is homotopy
equivalent to $\Omega^{2}(\mathcal{K}(\mathbb{R}))$, while $_{-1}%
\mathcal{V}(\mathbb{C})$ is homotopy equivalent to $\Omega^{3}(\mathcal{K}%
(\mathbb{R}))$.\hfill$\Box\smallskip$
\end{lemma}

Since all maps are module maps, our geometric viewpoint shows us that the
required morphism $\Omega^{2}(\mathcal{K}(\mathbb{R}))\longrightarrow
\Omega^{3}(\mathcal{K}(\mathbb{R}))$ is induced by the cup-product with an
element in the group $K_{1}(\mathbb{R})=\mathbb{Z}/2$. The following lemma
resolves this ambiguity and describes the homotopy fiber of the morphism.

\begin{lemma}
\label{_-1V(R) -> _-1V(C)}The map%
\[
{}_{-1}\mathcal{V}(\mathbb{R})\simeq\Omega^{2}(\mathcal{K}(\mathbb{R}%
))\longrightarrow{}_{-1}\mathcal{V}(\mathbb{C})\simeq\Omega^{3}(\mathcal{K}%
(\mathbb{R}))
\]
is induced by the cup-product with the generator of $K_{1}(\mathbb{R})$ and
its fiber has the homotopy type of $\mathcal{K}(\mathbb{C}).$ Therefore, the
homotopy fiber ${}_{-1}\mathcal{G}$ of the induced map on connective covers%
\[
{}_{-1}\mathcal{V}(\mathbb{R})^{c}\simeq(\Omega^{2}(\mathcal{K}(\mathbb{R}%
))^{c}\longrightarrow{}_{-1}\mathcal{V}(\mathbb{C})^{c}\simeq(\Omega
^{3}(\mathcal{K}(\mathbb{R})))^{c}%
\]
has the homotopy type of $\mathcal{K}(\mathbb{C})^{c}.$
\end{lemma}

\noindent\textbf{Proof.} To decide which element of $K_{1}(\mathbb{R})$ is
involved, one may use the fundamental theorem of hermitian $K$-theory in a
topological context (which is equivalent to Bott periodicity; see
\cite{Karoubi LNM343}). In other words, we can work in ${}_{1}U$-theory
instead of ${}_{-1}V$-theory. More precisely, if we show that the map ${}%
_{1}U_{0}(\mathbb{R})\longrightarrow{}_{1}U_{0}(\mathbb{C})$ is nontrivial,
this implies that our original map ${}_{-1}\mathcal{V}(\mathbb{R}%
)\rightarrow{}_{-1}\mathcal{V}(\mathbb{C})$ is \textit{not} nullhomotopic and
is therefore defined by the cup-product with the nontrivial element in
$K_{1}(\mathbb{R})$. For this purpose, we form the diagram%
\[%
\begin{array}
[c]{ccc}%
K_{1}(\mathbb{R}) & \rightarrow{}_{1}KQ_{1}(\mathbb{R}) & \rightarrow{}%
_{1}U_{0}(\mathbb{R})\\
\downarrow & \downarrow & \downarrow\\
K_{1}(\mathbb{C}) & \rightarrow{}_{1}KQ_{1}(\mathbb{C}) & \rightarrow{}%
_{1}U_{0}(\mathbb{C})
\end{array}
\]
with exact rows. With the aid of Lemmas \ref{_1KQ(R) splits} and
\ref{_1KQ(C)}, this can be rewritten as%
\[%
\begin{array}
[c]{cccc}%
\mathbb{Z}/2 & \rightarrow & \mathbb{Z}/2\oplus\mathbb{Z}/2 & \rightarrow
{}_{1}U_{0}(\mathbb{R})\\
\downarrow &  & \downarrow & \downarrow\\
0 & \rightarrow & \mathbb{Z}/2 & \rightarrow{}_{1}U_{0}(\mathbb{C})
\end{array}
\]
It is important to notice that the map%
\[
{}_{1}\mathcal{KQ}(\mathbb{R})\simeq\mathcal{K}(\mathbb{R})\times
\mathcal{K}(\mathbb{R})\longrightarrow{}_{1}\mathcal{KQ}(\mathbb{C}%
)\simeq\mathcal{K}(\mathbb{R})
\]
is the sum map. Therefore, the map $\mathbb{Z}/2\oplus\mathbb{Z}%
/2\rightarrow\mathbb{Z}/2$ in the diagram above is surjective. Since by
exactness ${}_{1}KQ_{1}(\mathbb{C})\rightarrow{}_{1}U_{0}(\mathbb{C})$ is
injective, it follows that {}$_{1}KQ_{1}(\mathbb{R})\rightarrow{}_{1}%
U_{0}(\mathbb{C})$ is nontrivial, and thus {}$_{1}U_{0}(\mathbb{R}%
)\rightarrow{}_{1}U_{0}(\mathbb{C})$ is nontrivial too.

Thus we obtain a nontrivial fiber sequence ${}_{-1}\mathcal{V}(\mathbb{R}%
)\longrightarrow{}_{-1}\mathcal{V}(\mathbb{C})$ with homotopy fiber
$\Omega^{2}(\mathcal{K}(\mathbb{C}))\simeq\mathcal{K}(\mathbb{C})$ by taking a
double loop of the Bott fibration%
\[
\mathcal{K}(\mathbb{C})\longrightarrow\mathcal{K}(\mathbb{R})\longrightarrow
\Omega(\mathcal{K}(\mathbb{R}))
\]
whose last map is defined by the cup-product with the nontrivial element in
$K_{1}(\mathbb{R}).$ Therefore, we also have a fibration of connective covers%
\[
\mathcal{K}(\mathbb{C})^{c}\longrightarrow{}_{-1}\mathcal{V}(\mathbb{R}%
)^{c}\longrightarrow{}_{-1}\mathcal{V}(\mathbb{C})^{c}%
\]
since ${}_{-1}V_{0}(\mathbb{C})=0$ by the previous lemma.\hfill$\Box
\smallskip$

\section{A ring with trivial $K$-theory and nontrivial $KQ$-theory}

Our first purpose here is to construct, from any ring $A$ in which $2$ is
invertible, a ring $R_{\infty}=R_{\infty}(A)$ with the following properties:

\begin{enumerate}
\item[(i)] \emph{For all }$n\in\mathbb{Z}$\emph{, }$K_{n}(R_{\infty}%
)=0$\emph{; but }

\item[(ii)] \emph{not all groups }${}_{\varepsilon}KQ_{n}(R_{\infty})$\emph{
need be trivial. }\newline\emph{In particular, for }$A$\emph{ a field of
characteristic }$\neq2$ \emph{and }$\varepsilon=1$\emph{,} $\emph{we}$
\emph{have }${}_{\varepsilon}KQ_{0}(R_{\infty})\cong W(A)$\emph{, the Witt
group of }$A$\emph{.}
\end{enumerate}

The existence of such a ring $R_{\infty}$ is used below to provide a
counterexample to a conjecture of \cite[3.4.2]{Williams}.

First recall that the suspension $S\Lambda$ of a discrete ring $\Lambda$ is
defined to be the quotient ring $C\Lambda/\tilde{\Lambda}$, where the cone
$C\Lambda$ of $\Lambda$ is the ring of infinite matrices (indexed by
$\mathbb{N}$) over $\Lambda$ for which there exists a natural number that bounds:

\begin{enumerate}
\item[(i)] the number of nonzero entries in each row and column; and

\item[(ii)] the number of distinct entries in the entire matrix.
\end{enumerate}

\noindent The ideal $\tilde{\Lambda}$ of $C\Lambda$ comprises matrices with
only finitely many nonzero entries.

Writing $\mathbb{Z}^{\prime}=\mathbb{Z}[1/2]$, for a $\mathbb{Z}^{\prime}%
$-algebra $A$ denote by $\mathcal{P}(A)$ the category of its finitely
generated projective right modules. The tensor product of modules over
$\mathbb{Z}^{\prime}$ then defines a biadditive functor
\[
\mathcal{P}(A)\times\mathcal{P}(S^{2}\mathbb{Z}^{\prime})\longrightarrow
\mathcal{P}(S^{2}A)\text{,}%
\]
where $S^{2}$ refers to the double suspension.

Now provide $M_{2}(S^{2}\mathbb{Z}^{\prime})$ and $M_{2}(S^{2}A)$ with the
involution%
\[
\left[
\begin{tabular}
[c]{cc}%
$a$ & $b$\\
$c$ & $d$%
\end{tabular}
\right]  \longmapsto\left[
\begin{tabular}
[c]{cc}%
$\overline{d}$ & $-\overline{b}$\\
$-\overline{c}$ & $\overline{a}$%
\end{tabular}
\right]
\]
and choose a self-adjoint projection operator $p$ in $M_{2}(S^{2}%
\mathbb{Z}^{\prime})$ whose image defines a class in $_{1}KQ_{0}(M_{2}%
(S^{2}\mathbb{Z}^{\prime}))\cong\mathsf{{}}_{-1}KQ_{-2}(\mathbb{Z}^{\prime
})\cong\mathbb{Z}\oplus\mathbb{Z}/2$ that is a generator of the free summand
\cite{K:AnnM112hgo}. The tensor product with $p$ induces a nonunital map
between unital rings%
\[
\phi:A\longrightarrow M_{2}(S^{2}A)
\]
defined by $a\mapsto a\otimes p$. It is easy to see that $\phi$ induces on
$K_{0}$ and $KQ_{0}$ the cup-product with the class in $_{-1}KQ_{-2}%
(\mathbb{Z}^{\prime})$ mentioned above.

In order to deal with the technical problem that $\phi(1)\neq1$, let us
replace the ring $A$ by the suspension $SB$ of a ring $B$. We can view $A$ as
a bimodule over the cone $CB$ of this ring $B$. This allows us to
\textquotedblleft add a unit\textquotedblright\ to the ring $A$ by defining
$R$ as $CB\times A$ with the multiplication rule%
\[
(\lambda,u)(\mu,v)=(\lambda\mu,\,\lambda v+u\mu+uv)\text{.}%
\]
We obtain an exact sequence of rings and nonunital ring homomorphisms%
\[
0\longrightarrow A\longrightarrow R\longrightarrow CB\longrightarrow0
\]
which shows that in negative degrees the $K$-theories of $R$ and $A$ coincide,
as do their $KQ$-theories. We then change $\phi$ into a map between unital
rings%
\[
\Phi:R\longrightarrow R\otimes_{\mathbb{Z}^{\prime}}S^{2}(\mathbb{Z}^{\prime
})=S^{2}(R)=R_{1}%
\]
by the formula $\Phi(\lambda,u)=(\lambda\otimes1,\,u\otimes p)$. This map
restricts to $\phi$ on the ideal $A$ and we can safely use $\Phi$ as a
substitute for $\phi$. Thus, we are able to define a direct system of unital
rings with involution $R_{t}$ by the inductive formula $R_{t+1}=(R_{t})_{1}$.

By a well-known theorem of Wagoner \cite{Wagoner} (see also
\cite{K:AnnM112fun}), we have canonical isomorphisms $K_{n+1}(SD)\cong
K_{n}(D)$ and ${}_{\varepsilon}KQ_{n+1}(SD)\cong{}_{\varepsilon}KQ_{n}(D)$.
Therefore, if we define $R_{\infty}$ as the direct limit of the $R_{t}$, we
have%
\[
K_{n}(R_{\infty})\cong\underrightarrow{\lim}{}K_{n}(R_{t})\cong%
\underrightarrow{\lim}{}K_{n-2t}(A)=0
\]
since $K_{-2}(\mathbb{Z}^{\prime})=0$.\textsf{ }

On the other hand ${}_{\varepsilon}KQ_{n}(R_{\infty})\cong\underrightarrow
{\lim}{}{}_{\varepsilon}KQ_{n-4s}(A)$ is the stabilized Witt group of $A$
\cite{Karoubi stab.Witt}, which is not trivial in general. For instance, if
$A$ is a commutative regular noetherian ring and $\varepsilon=1$, this is the
classical Witt group of $A$. Hence, our construction of the ring $R_{\infty}$
is complete.

This construction provides many counterexamples to a conjecture of B. Williams
\cite{Williams}. More specifically, we have for instance the following theorem.

\begin{theorem}
Let $A$ be a commutative regular noetherian ring with finitely generated
$\varepsilon$-Witt groups in degrees $0$ and $1$, and let $R_{\infty}$ be the
associated ring defined above. Then the canonical map%
\[
{}_{1}\mathcal{KQ}(R_{\infty})_{\#}^{c}\longrightarrow(\mathcal{K}(R_{\infty
})_{\#}{}^{h({}_{1}\mathbb{Z}/2)})^{c}%
\]
is NOT a homotopy equivalence.
\end{theorem}

\noindent\textbf{Proof.} According to the computation before (using Lemma
\ref{Connective.Spectra}\thinspace(v)), $\pi_{0}({}_{1}\mathcal{KQ}(R_{\infty
})_{\#}^{c})$ is the $2$-completed Witt group $W(A)_{\#}$ because
$W(R_{\infty})=W(A)$ in this case \cite{Karoubi stab.Witt}. The group
$W(A)_{\#}$ is not trivial since the rank map induces a surjection from this
group to $\mathbb{Z}/2$. On the other hand, since $\mathcal{K}(R)_{\#}$ has
trivial homotopy groups, it is contractible, which implies that the group
$\pi_{0}$ of the right hand side is reduced to $0$.\hfill$\Box\smallskip$

\section{Adams operations on the real $K$-theory spectrum}

Here we consider the composite%
\[
\Omega^{2}\mathcal{K}(\mathbb{R})\overset{\sigma}{\longrightarrow}\Omega
^{3}\mathcal{K}(\mathbb{R})\overset{\Omega^{3}(\psi^{q}-1)}{\longrightarrow
}\Omega^{3}\mathcal{K}(\mathbb{R})
\]
where $q$ is odd and $\sigma$ is induced by the cup-product with the generator
$H-1$ of $K_{1}(\mathbb{R})=\mathbb{Z}/2$, and the canonical real line bundle
$H$ over $S^{1}$ has $H^{2}=1$.

\begin{proposition}
The composite%
\[
\Omega^{2}\mathcal{K}(\mathbb{R})_{\#}^{c}\overset{\sigma}{\longrightarrow
}(\Omega^{3}\mathcal{K}(\mathbb{R}))_{\#}^{c}\overset{\Omega^{3}(\psi^{q}%
-1)}{\longrightarrow}(\Omega^{3}\mathcal{K}(\mathbb{R}))_{\#}^{c}%
\]
is NOT nullhomotopic.
\end{proposition}

\noindent\textbf{Proof.\ }We first show that for any space $X$, the
composition of homotopy groups (where $\left[  X,\mathcal{E}\right]  $ refers
to pointed homotopy classes from $X$ to $\mathcal{E}_{0}$)%
\[
\left[  X,\,\mathcal{K}(\mathbb{R})\right]  \overset{\sigma_{\ast}%
}{\longrightarrow}\left[  X,\,\Omega\mathcal{K}(\mathbb{R})\right]
\overset{\Omega(\psi^{q}-1)_{\ast}}{\longrightarrow}\left[  X,\,\Omega
\mathcal{K}(\mathbb{R})\right]
\]
is equal to the composite in reverse order%
\[
\left[  X,\,\mathcal{K}(\mathbb{R})\right]  \overset{(\psi^{q}-1)_{\ast}%
}{\longrightarrow}\left[  X,\,\mathcal{K}(\mathbb{R})\right]  \overset
{\sigma_{\ast}}{\longrightarrow}\left[  X,\,\Omega\mathcal{K}(\mathbb{R}%
)\right]
\]
(in other words, $\sigma$ commutes with $\psi^{q}$). This follows from the
following straightforward computation:%
\[
\psi^{q}(\sigma(x))=\psi^{q}((H-1)\cdot x)=\psi^{q}(H-1)\cdot\psi
^{q}(x)=(H^{q}-1)\cdot\psi^{q}(x)=(H-1)\cdot\psi^{q}(x)=\sigma(\psi^{q}(x)),
\]
where $x\in K_{\mathbb{R}}(X).$

Therefore, on passing to $2$-completions of connective covers, our claim will
follow from the nontriviality of the composition%
\[
(\Omega^{2}\mathcal{K}(\mathbb{R}))_{\#}^{c}\overset{\Omega^{2}(\psi^{q}%
-1)}{\longrightarrow}(\Omega^{2}\mathcal{K}(\mathbb{R}))_{\#}^{c}%
\overset{\sigma}{\longrightarrow}(\Omega^{3}\mathcal{K}(\mathbb{R}))_{\#}^{c}%
\]
or%
\begin{equation}
(\Omega^{8}\mathcal{K}(\mathbb{R}))_{\#}^{c}\overset{\Omega^{8}(\psi^{q}%
-1)}{\longrightarrow}(\Omega^{8}\mathcal{K}(\mathbb{R}))_{\#}^{c}%
\overset{\sigma}{\longrightarrow}\Omega^{9}\mathcal{K}(\mathbb{R})_{\#}^{c}.
\label{Adams}%
\end{equation}

By a theorem of Bott (see for instance \cite[Section III.5]{K:ktheorybook}),
there is a fiber sequence%
\[
\mathcal{K}(\mathbb{C})\overset{r}{\longrightarrow}\mathcal{K}(\mathbb{R}%
)\overset{\sigma}{\longrightarrow}\Omega\mathcal{K}(\mathbb{R})
\]
where the homotopy fiber of $\sigma$ is $\mathcal{K}(\mathbb{C})$, the
classifying space of complex topological $K$-theory, and $r$ is the
realification map. Therefore, if the sequence (\ref{Adams}) were trivial, one
would have a factorization%
\[
\Omega^{8}(\psi^{q}-1):(\Omega^{8}\mathcal{K}(\mathbb{R}))_{\#}^{c}%
\longrightarrow(\Omega^{8}\mathcal{K}(\mathbb{C}))_{\#}^{c}\overset
{r}{\longrightarrow}(\Omega^{8}\mathcal{K}(\mathbb{R))}_{\#}^{c}%
\]
or equally, by Bott periodicity,%
\[
q^{4}\psi^{q}-1:\mathcal{K}(\mathbb{R})_{\#}^{c}\longrightarrow\mathcal{K}%
(\mathbb{C})_{\#}^{c}\overset{r}{\longrightarrow}\mathcal{K}(\mathbb{R}%
)_{\#}^{c}.
\]

A way to prove the impossibility of such a factorization is to find a test
space $X$ and map it into the three spaces involved. For such a space we
choose the classifying space $X=$ $BG$, where $G$ is the connected Lie group
$SO(3)$. We thereby transform our problem into an algebraic one: by Atiyah and
Hirzebruch \cite[Theorem 4.8]{Atiyah-Hirzebruch}, we know that $K_{\mathbb{C}%
}(BG)=\widehat{R(G)}$, the complex representation ring of $G$ completed at its
augmentation ideal, while by Anderson \cite{Anderson}, we have $K_{\mathbb{R}%
}(BG)=\widehat{RO(G)}$, the completed real representation ring of $G$.

It is well-known (for example, \cite{Adams}) that $RO(G)$ is the polynomial
algebra in one variable $\mathbb{Z}\left[  \lambda^{1}\right]  $, where
$\lambda^{1}$ is the standard representation of $SO(3)$ in $\mathbb{R}^{3}$.
The complexification $c$ induces an isomorphism $RO(G)\overset{\cong%
}{\rightarrow}R(G)$, while the realification from $R(G)$ to $RO(G)\cong R(G)$
is identified with the multiplication by $2$ (any complex representation of
$G$ is isomorphic to its conjugate). If we choose $SO(2)$ as a maximal compact
torus in $SO(3)$ embedded as%
\[
\left[
\begin{array}
[c]{ccc}%
\cos\theta & -\sin\theta & 0\\
\sin\theta & \cos\theta & 0\\
0 & 0 & 1
\end{array}
\right]
\]
we may view $R(G)$ as the ring of polynomials on the variable $\lambda
^{1}=t+t^{-1}+1$, where $t$ represents the standard one-dimensional
representation $\theta\mapsto e^{i\theta}$ of $SO(2)=S^{1}$.

Now, since $G$ is a compact connected Lie group, $R(G)$ injects into its
completion. If we put $t=1-u$, we may identify the completed representation
ring $\widehat{R(G)}$ as a subring of the ring of formal power series
$\mathbb{Z}\left[  \left[  u\right]  \right]  $. However, we are considering
homotopy classes from $X=BG$ not just to $\mathcal{K}(\mathbb{R})_{0}$ or
$\mathcal{K}(\mathbb{C})_{0}$ but to its $2$-adic completion $\mathcal{K}%
(\mathbb{R})_{0\#}$ or $\mathcal{K}(\mathbb{C})_{0\#}$ \cite[(v),
pg.~205]{Adamsbluebook}. From the algebraic point of view, this means that we
have to compute in the power series ring $\mathbb{Z}_{2}\left[  \left[
u\right]  \right]  $ instead of $\mathbb{Z}\left[  \left[  u\right]  \right]
$.

Since $t$ is one-dimensional and $\psi^{q}$ commutes with the complexification
isomorphism $c$ between $R(G)$ and $RO(G)$, we can write, in $\mathbb{Z}%
_{2}\left[  \left[  u\right]  \right]  $:%
\begin{align*}
c(q^{4}\psi^{q}-1)(\lambda^{1}) &  =q^{4}t^{q}+q^{4}t^{-q}+q^{4}-t-t^{-1}-1\\
&  =(1-u)^{-q}[q^{4}(1-u)^{2q}-(1-u)^{q+1}+(q^{4}-1)(1-u)^{q}-(1-u)^{q-1}%
+q^{4}]\text{.}%
\end{align*}
Because multiplication by $(1-u)^{q}$ leaves an odd coefficient of $u^{2q}$,
the power series is not divisible by $2$. Therefore $(q^{4}\psi^{q}%
-1)(\lambda^{1})$ cannot be in the image of the realification map $r$ (which
furnishes elements divisible by $2$ in $\widehat{RO(G)}\overset{\underset
{\cong}{c}}{\longrightarrow}\widehat{R(G)}$). This contradiction concludes the
proof of the proposition.$\hfill\Box\smallskip$

\bigskip

\section*{Acknowledgements}

For hospitality and support during the development of this paper, the first
author thanks the Departments of Mathematics at the Universities Paris 7 and
Blaise Pascal Clermont-Ferrand as well as NUS research grant
R-146-000-097-112. The second and third authors gratefully acknowledge the
Thematic Program in Homotopy Theory at the Fields Institute, and the third
author is pleased to extend his thanks to the Homotopy Theory and Higher
Categories programme organized by Centre de Recerca Matematica at Universitat
Autonoma de Barcelona.

We would also like to thank the referees of our paper very warmly. Their
detailed comments have contributed, we believe, to a better presentation of
the results and proofs which are now more precise and detailed. Thanks to this
second careful reading, we corrected Theorem \ref{theorem4} from the first
version. This has led to the contents of Appendix D dealing with the
symplectic case becoming more interesting. As a drawback, our paper is longer
than in the first version found in the ArXiv, but we hope more readable for
the majority of the readers.



\begin{center}
A.\thinspace Jon Berrick \\[0pt]Department of Mathematics, National University
of Singapore, Singapore. \\[0pt]e-mail: berrick@math.nus.edu.sg \vspace{0.2in}

Max Karoubi \\[0pt]UFR de Math{\'{e}}matiques, Universit{\'{e}} Paris 7,
France. \\[0pt]e-mail: max.karoubi@gmail.com \vspace{0.2in}

Paul Arne {\O }stv{\ae }r \\[0pt]Department of Mathematics, University of
Oslo, Norway. \\[0pt]e-mail: paularne@math.uio.no
\end{center}

\end{document}